\documentstyle[12pt]{article}
\begin{document}
\input{amssym.def}
\input{amssym}


\def\theequation{\thesection.\arabic{equation}}
\def\be{\begin{equation}}
\def\ee{\end{equation}}
\def\ba{\begin{eqnarray}}
\def\ea{\end{eqnarray}}
\def\lb{\label}
\def\nn{\nonumber}

\def\a{\alpha}
\def\b{\beta}
\def\d{\delta}
\def\e{\varepsilon}
\def\E{{\cal E}}
\def\bq{\overline{q}}

\def\id{\mbox{\rm 1\hspace{-3.5pt}I}}
\newcommand{\ID}[2]{\id^{| #1 {\cal i}}_{\;\;\; {\cal h} #2 |}}

\def\C{\Bbb C}
\def\Z{\Bbb Z}
\def\F{\Bbb F}
\def\eod{\phantom{a}\hfill \rule{2.5mm}{2.5mm}}

\def\hR{\hat{R}}\def\hF{\hat{F}}
\def\hA{\hat{A}}\def\hB{\hat{B}}
\newcommand{\dy}[1]{DYBE(${#1}$)}
\newcommand{\xR}[1]{\ ^{\small #1}\!\!\hR}
\newcommand{\xF}[1]{\ ^{\small #1}\!\!\hF}

\newcommand{\rank}{\mathop{\rm rank}\nolimits}
\newcommand{\height}{\mathop{\rm height}\nolimits}
\newcommand{\aut}{\mathop{\rm Aut}\nolimits}
\newcommand{\rx}{\mathop{\rho_{\hspace{-1pt}\scriptscriptstyle W,k}}\nolimits}
\newcommand{\rn}{\mathop{\rho_{\hspace{-1pt}\scriptscriptstyle W,n}}\nolimits}
\newcommand{\q}[1]{[#1]}
\newcommand{\ai}[1]{{a_{#1}}\, }
\newcommand{\ainv}[2]{({a^{-1})}^{| #1 {\cal i}}_{\;\;\; {\cal h} #2 |}}

\def\R{\hat{R}}
\def\Rp{\hat{R}(p)}
\newcommand{\DR}[1]{\hat{R}_{#1}(p)}
\newcommand{\DDR}[2]{\hat{R}^{#1}_{#2}(p)}

\def\A{{A}}
\def\eup{\varepsilon^{|1  \dots n{\cal i}}}
\def\edo{\varepsilon_{{\cal h}1  \dots n|}}
\def\eu2{\varepsilon^{|2  \dots n{+}1{\cal i}}}
\def\ed2{\varepsilon_{{\cal h}2  \dots n{+}1|}}
\def\Eup{{\cal E}^{|1  \dots n{\cal i}}(p)}
\def\Edo{{\cal E}_{{\cal h}1  \dots n|}(p)}
\def\Eu2{{\cal E}^{|2  \dots n{+}1{\cal i}}(p)}
\def\Ed2{{\cal E}_{{\cal h}2  \dots n{+}1|}(p)}

\newcommand{\eupi}[1]{\varepsilon^{|#1 \dots n+ #1 -1{\cal i}}}
\newcommand{\edoi}[1]{\varepsilon_{{\cal h}#1 \dots n+ #1 -1|}}
\newcommand{\Eupi}[1]{{\cal E}^{|#1 \dots n+ #1 -1{\cal i}}(p)}
\newcommand{\Edoi}[1]{{\cal E}_{{\cal h}#1 \dots n+ #1 -1|}(p)}
\newcommand{\N}[2]{N^{| #1 {\cal i}}_{\;\;\; {\cal h} #2 |}}
\newcommand{\K}[2]{K^{| #1 {\cal i}}_{\;\;\;\;\;\;\;\: {\cal h} #2 |}}
\newcommand{\iK}[2]{{K^{-1}}^{| #1 {\cal i}}_{\;\;\; {\cal h} #2 |}}

\def\bbr{{\rm I}\!{\rm R}}
\def\bbz{{\rm Z}\!\!\!{\rm Z}}
\def\subbbc{{\rm C}\kern-3.3pt
\hbox{\vrule height4.8pt width0.4pt}\, }
\def\BbbZ{Z\!\!\!Z}
\def\BbbN{{\rm I}\!{\rm N}}

\hyphenation{dy-na-mi-cal}
\hyphenation{ge-ne-ra-ted}

\begin{center}


{\Large Hecke algebraic properties of}\\[3 mm]
{\Large dynamical $R$-matrices. Application}\\[3 mm]
{\Large to related quantum matrix algebras}\\[8 mm]

{\large L.K. Hadjiivanov$^{a}$
\footnote[1]
{On leave of absence from:
Division of Theoretical Physics, Institute for Nuclear
Research and Nuclear Energy,
Bulgarian Academy of Sciences, Tsarigradsko Chaussee
72, BG-1784 Sofia, Bulgaria; e-mail address: lhadji@inrne.acad.bg},
A.P. Isaev$^{b}$
\footnote[2]
{On leave of absence from: Bogoliubov Laboratory of Theoretical Physics, JINR,
Dubna, 141 980 Moscow Region, Russia;  e-mail address:
isaevap@thsun1.jinr.ru},
O.V. Ogievetsky$^{c}$
\footnote[3]
{
On leave of absence from: P.N. Lebedev Physical Institute,
Theoretical Department, 117924 Moscow, Leninsky prospect 53, Russia;
e-mail address: oleg@cpt.univ-mrs.fr},\\[.2 cm]
P.N. Pyatov$^{d}$
\footnote[4]
{e-mail address: pyatov@thsun1.jinr.ru}
{\normalsize and}
I.T. Todorov$^{e}$
\footnote[5]
{On leave
of absence from: Division of Theoretical Physics, Institute for Nuclear
Research and
Nuclear Energy, Bulgarian Academy of Sciences, Tsarigradsko Chaussee
72, BG-1784 Sofia, Bulgaria; e-mail address: todorov@inrne.acad.bg} }\\[.6 cm]
{\footnotesize $^{a}$International
Centre for Theoretical Physics (ICTP), I-34014
Trieste, Italy }\\
{\footnotesize $^{b}$Dipartimento di Fisica, Universit\`a di Pisa,
I-56100 Pisa, Italy }\\
{\footnotesize $^{c}$Centre de Physique Th\'eorique, Luminy, F-13288
Marseille, France }\\
{\footnotesize $^{d}$Bogoliubov Laboratory of Theoretical Physics, JINR,}\\
{\footnotesize Dubna, 141980 Moscow Region, Russia }\\
{\footnotesize $^{e}$Erwin Schr\"odinger Institute for Mathematical
Physics (ESI),}\\
{\footnotesize A-1090 Wien, Austria }\\[.2 cm]
\end{center}


\begin{abstract}
{\normalsize
\noindent
The quantum dynamical Yang--Baxter (or Gervais--Neveu--Felder)
equation defines an $R$-matrix $\R(p)\,$, where $p$ stands for a
set of mutually commuting variables.  A family of $SL(n)$-type
solutions of this equation provides a new realization of the Hecke
algebra. We define quantum antisymmetrizers, introduce the notion of
quantum determinant and compute the inverse quantum matrix for matrix
algebras of the type
${\hat R}(p) a_1 a_2 = a_1 a_2 {\hat R}$.
It is pointed out that such a quantum matrix algebra arises in
the operator realization of the chiral zero modes of the WZNW model.
}
\end{abstract}

\newpage


\textwidth = 16truecm
\textheight = 24truecm
\hoffset = -1truecm
\voffset = -2truecm

\section*{Introduction}
\setcounter{equation}{0}
\renewcommand\theequation{0.\arabic{equation}}

\medskip

Let $\{\,v^{(i)}\,,\,i=1,\ldots ,n\,\}$ be a "barycentric basis"
in a Cartan subalgebra ${\frak h} \subset sl(n)\,.\,$ Viewed as operators
in the $n$-dimensional complex space $V = {\C}^n\,,\, v^{(i)}\,$ can be
realized as real traceless diagonal $n\times n$ matrices:
\be
\lb{0.1}
(v^{(i)})^j_j = \d_{ij} - {1\over n}\,,\quad
\sum_{i=1}^n\, v^{(i)} = 0\,.
\ee
Let further $\{ p_i\}_{i=1}^n\,$ span the dual Lie algebra ${\frak h}^*\,.$
Introduce the traceless diagonal matrix
\be
\lb{0.2}
p = p_i v^{(i)} \ \left(\equiv \sum_{i=1}^n p_i v^{(i)} \right)\,,\quad
[p_i\,,\, p_j\, ]=0\,, \quad \sum_{i=1}^n p_i = 0\,.
\ee
We define a Hecke-type quantum dynamical $R$-matrix
$\Rp$ as a map from ${\frak h}^*$ to $End\;\left( V\otimes V\right)$
satisfying the {\em twisted braid relation}
\be
\lb{0.3}
\DR{12} \R_{23}(p-v_1) \DR{12} = \R_{23}(p-v_1)\DR{12} \R_{23}(p-v_1)
\ee
and the Hecke condition
\be
\lb{0.4}
\Rp^2 = \id + (q- {\bar q} )\Rp\,,\quad {\bar q}\, := \, q^{-1}\,.
\ee
(Although the notation is taylored to the special case in which
the parameter $q\,$ takes values on the unit circle, we shall not use this
property in the main body of the paper.) The subscripts in (\ref{0.3})
refer to the, by now standard, tensor product notation of Faddeev et al.
(see, e.g., \cite{FRT}); in particular, $\R_{23}(p-v_1)\, \in\,
End\;(V^{\otimes 3} )$
has matrix elements
\be
\lb{0.5}
\left(\R_{23}(p-v_1) \right)^{i_1 i_2 i_3}_{j_1 j_2 j_3}\, =\,
\d_{j_1}^{i_1}\; \R (p-v^{(i_1)})^{i_2 i_3}_{j_2 j_3}\,.
\ee
The twisted braid relation (\ref{0.3}) is equivalent to the
{\em quantum dynamical} \cite{GN,Skl} (or {\em deformed} \cite{BF})
{\em Yang--Baxter equation} (QDYBE) for the matrix $R(p)$ related to the braid
operator $\Rp$ by $\Rp = P R(p)$ where $P$ stands for the permutation
operator
$P x_1 y_2  = y_1 x_2 ,\ x, y \in  V ,\  P^2 = 1$.
Abusing notation we
shall also refer to Eq.(\ref{0.3}) by the above abbreviation.
The term "dynamical $R$-matrix" for $\Rp$ is suggested by the fact that in
the physical applications its arguments play the role of (commuting)
dynamical variables and that $\Rp$ satisfies a finite difference rather
than a purely algebraic equation.

The important concept of a {\em quantum matrix algebra}
${\cal A} = {\cal A}(\Rp, \R)$ (Sec.5) can be
introduced as a (complex) associative algebra with $1$
generated by rational functions of $q^{p_i},~i=1,\dots ,n$, and the (noncommuting) entries of an
$n\times n$ matrix $a = \left( a^i_\alpha
\right)\,$ satisfying the quadratic exchange relations
\be
\lb{0.6}
\Rp \, a_1 \, a_2  \, =  \, a_1 \, a_2 \, \R \; ,
\ee
where $\R \equiv \R^{\a_1 \a_2}_{\b_1 \b_2}$
in the right hand side is a constant (i.e., $p$-independent)
solution of (\ref{0.3}), (\ref{0.4}), and all entries in a matrix row
$a^i\,= \{a^i_\a\}_{\a=1}^n$ are acting equivalently as
shift operators for $p\,$:
\be
\lb{0.7}
p\, a^i\, =\, a^i\, (\, p \,+\, v^{(i)}\, )\ \  {\rm or} \ \
p_{jk}\, a^i\, = \, a^i\, (p_{jk} + \d^i_j - \d^i_k \,)\ \ {\rm for}\ \
p_{jk} = p_j - p_k\,.
\ee
\vspace{0.2 cm}

\noindent
{\bf Remark 1~}
It has been pointed out \cite{AF, BF} that, in the $su(2)$ case, the
matrix $a$ generates the $q$-Clebsh--Gordan coefficients while
$\Rp$ plays the role of a "quantum $6j$-symbol" \cite{F1,KR,BBB}.

\vspace{0.2 cm}


\noindent
{\bf Remark 2~}
Eq.(\ref{0.6}) is related to the one with indices $1$ and $2$ interchanged,
\be
\lb{0.10}
\Rp \, a_2 \, a_1 \, = \, a_2 \, a_1 \, \R \; ,
\ee
by the substitution $\R \rightarrow P \R P\,,\,\,\Rp \rightarrow P \Rp\,P\,$.
It is, on the other hand, formally obtained from
\be
\lb{0.11}
\R \, \overline{a}_2 \, \overline{a}_1 \, = \,
\overline{a}_2 \, \overline{a}_1 \, \Rp \; .
\ee
by the substitution $\overline{a} = a^{-1}\,$; the same
substitution relates (\ref{0.10}) to
\be
\lb{0.12}
\R \, \overline{a}_1 \, \overline{a}_2 \, = \,
\overline{a}_1 \, \overline{a}_2 \, \Rp \;.
\ee
Since
\be
\lb{PRP}
{\Rp}_{21} = P\,\Rp P
\ee
satisfies conditions of the same type as $\Rp\,,$
we can start with either
of these relations.
\vspace{0.2 cm}

The QDYBE, introduced by Gervais and Neveu \cite{GN} for the exchange algebra
associated with the Liouville equation and applied to the
zero mode algebra \cite{AF} of the Wess--Zumino--Novikov--Witten
(WZNW) model \cite{W,F11,F2} is attracting
ever more attention.
Its classical
counterpart, introduced in \cite{Skl} (see also \cite{Mai}) has been
displayed in \cite{AF} for the $sl(2)$ case and in \cite{AT} for an arbitrary
simple Lie algebra.
The quantum $\Rp$ is central to a continuing study \cite{BF} of
$q$-deformed cotangent bundles on group manifolds and quantum model spaces. A
quasiHopf-algebraic point of view is taken in \cite{By, BS} where $\Rp$ is
obtained by a Drinfeld twist of the constant $R$-matrix. Felder \cite{Fe}
explores the more general case of (classical and) QDYBE depending on a spectral
parameter and finds elliptic solutions of this equation. These solutions
are applied in \cite{ABB} to
quantize Calogero--Moser and Ruijsenaars--Schneider
models. A class of $SL(n\vert m )$-type solutions of the QDYBE (and
related trigonometric solutions of the equation with spectral parameter) are
described in \cite{I2}. A more mathematically minded
approach to the subject
in terms of "$\frak h$-algebroids" is being developed in \cite{EV}.

The present work was motivated in part by earlier study \cite{FHT1, FHT2,
DT} of the canonical quantization of the WZNW model (following \cite{F1,
F11, AF, F2, G}).
It was noticed, in particular, that the exchange relations
\cite{FHT2} for the chiral zero modes $a^i_\alpha\,$ that diagonalize the
$U_q(sl(2))\,$ monodromy matrix can be written in the form (\ref{0.10}).
As a result, the operator realization of the chiral group valued field was
understood as a quantization of the (deformed) classical Poisson bracket
relations of \cite{AT} thus opening the way to its generalization for
$SU(n)$. Here we show that
a special solution of the QDYBE (\ref{0.3}) yields a new
matrix representation of the Hecke algebra. We
concentrate on a general study of the Hecke algebra
properties of this solution
and the ensuing properties of the quantum matrices
satisfying (\ref{0.6}) relegating applications to the WZNW model to a
subsequent publication \cite{FHIPT} which is highlighted
in Section 6.
A central result is the computation of the quantum determinant of
$a$ and the (based on it) evaluation of the
inverse quantum matrix.

The paper is organized as follows. We review and extend in Section 1
results of Gurevich \cite{Gur} on quantum
(anti)symmetrizers
and illustrate them in Section 2 on the known example of a constant $\R\;$.
We proceed in Section 3 to a study of a family of $SL(n)$-type dynamical
$R$-matrices and describe two types of symmetry
transformations for this family: the twist transformation
(a version of Drinfeld's twist for dynamical $R$-matrices)
and the canonical shifts.
In Section 4 we show that
these dynamical matrices provide a new realization
of the Hecke algebra. This
allows to define "dynamical" ($p$-dependent) analogs of quantum
antisymmetrizers, including the Levi--Civit\`a $\cal E$ tensor. In
Section 5 we study the quantum algebra $\cal A\,$, define the quantum
determinant ${\rm det}\,(a)$ and compute the inverse matrix $a^{-1}\,$. We
demonstrate that $\cal A$ provides a realization of a reflection equation
algebra which is interpreted as a quantum monodromy algebra in the WZNW
theory. An Appendix is devoted to deriving some useful identities
for the parameters determining the solution of  the QDYBE found in
\cite{I2} and to computing the
normalization of the dynamical Levi--Civit\`a tensor.


\newpage

\section{Hecke algebras and $q$-antisymmetrizers}
\setcounter{equation}{0}
\renewcommand\theequation{\thesection.\arabic{equation}}

\medskip

In this Section we collect some basic notions on Hecke algebras
and describe the $q$-antisymmetrizers technique
which is to be applied later on.
We follow closely the approach of D.~Gurevich \cite{Gur}.
\vspace{3mm}

In the present context by a
Hecke algebra ${\cal H}_k(q)$ we understand a $\C$-algebra
with generators $1, g_1, g_2, \dots , g_{k-1}$, a nonzero
parameter $q \in \C$, and defining relations
\ba
\lb{a1}
g_i \, g_{i+1} \, g_i = g_{i+1} \, g_i \, g_{i+1} && \mbox{for} \;\;
1\leq  i \leq k-2, \\
\lb{a2}
{g_i}^2 = 1 + (q-\bq)\, g_i && \mbox{for}\;\; 1\leq i \leq k-1, \\
\lb{a3}
g_i \, g_j = g_j \, g_i && \mbox{if} \;\; |i-j|\geq 2,
\ea
where $\bq := q^{-1}$.
\medskip

We shall consider the set of
idempotents $\A^{(j)} \in {\cal H}_k(q),~j = 1,\dots, k$,
associated with single column
Young diagrams containing $j$ nodes -- the so called $q$-antisymmetrizers.
Their inductive definition is given by
\be
\lb{a4}
\A^{(1)} = 1\, , \quad
\A^{(j)} = {1 \over \q{j}}
\A^{(j-1)}\left( q^{j-1} - \q{j-1} g_{j-1} \right) \A^{(j-1)}\, .
\ee
Here $\q{j} = (q^j - {\bq}^j)/(q-\bq)$
and we assume $\q{j}\neq 0$, for $j=2,\dots ,k$.
Note, that $A^{(k)}$ is a central idempotent in the algebra
${\cal H}_{k}(q)$.

Equivalently, one may write
\be
\lb{a5}
\A^{(j)} = {1 \over \q{j}}
\A^{(2,j-1)} \left( q^{j-1} - \q{j-1} g_1 \right)\A^{(2,j-1)} \, ,
\ee
where we have adopted the notation $\A^{(i,j)}$, $1\leq i\leq j$, for
the central idempotent  of the subalgebra
${\cal H}_{i,j}(q) \subset {\cal H}_j(q)$
generated by the subset
$1, g_i,\dots ,g_{j-1}$.
In particular,  $\A^{(1,j)} = \A^{(j)}$, $\A^{(j,j)}= 1$.

\vspace{2mm}
\noindent
{\bf Remark 1.1~}
All the subalgebras
\mbox{${\cal H}_{i,r+i}(q)\subset {\cal H}_k(q)$}, \mbox{$i=1,\dots ,k{-}r$},
are isomorphic by definition.
Moreover, they are related by inner ${\cal H}_k(q)$-auto\-mor\-phisms.
For example, the automorphism
\mbox{$\phi_i : {\cal H}_{i,r+i}(q) \to {\cal H}_{i+1,r+i+1}(q)$} is given by
\be
\lb{aut}
\phi_i(t) = g_i g_{i+1}\dots g_{r+i}\, t\,
(g_i g_{i+1}\dots g_{r+i})^{-1}\ ,
\quad \forall t\in {\cal H}_{i,r+i}(q)\ .
\ee

The term $q$-antisymmetrizer for the elements $\A^{(j)}$
is justified by the following properties:
\ba
\lb{a6}
(g_i + \bq) \A^{(j)} = \A^{(j)} (g_i +\bq) = 0 &&
\mbox{for} \;\; 1 \leq i \leq j-1 \, , \\
\lb{a7}
\A^{(j)} \A^{(i,l)} = \A^{(i,l)} \A^{(j)} = \A^{(j)} && \mbox{for} \; \;
1\leq i\leq l\leq j\, .
\ea

\vspace{2mm}
\noindent
{\bf Remark 1.2~}
Replacing $q$ by $(-q^{-1})$ in (\ref{a4}) leads to another sequence
of projectors, called symmetrizers. Abstractly, inside the Hecke algebra,
it is a matter of convention - which projectors one calls symmetrizers,
and which - antisymmetrizers. We use the common convention. However,
on the level of representations, when one can calculate the ranks of
projectors and see which sequence of projectors terminates, the
distinction between symmetrizers and antisymmetrizers
becomes
meaningful.

\vspace{2mm}
Consider a representation
$\rx :\, {\cal H}_k(q) \rightarrow \aut(W)$
of the algebra ${\cal H}_k(q)$ in a vector space $W$.

\vspace{2mm}
\noindent
{\bf Definition 1.1~}
{\em
We shall say that $\rx$	is a representation of $\height n$
in one of the following two cases:

{\bf a)} $n < k$ and the conditions
\ba
\lb{a8}
\rx(\A^{(n+1)}) &=& 0\ ,
\\
\lb{a9}
\rank \rx(\A^{(n)}) &=& 1\ ,
\ea

\vspace{2mm}

are fulfilled, or

\vspace{2mm}

{\bf b)} $n = k$ and \hspace{1.45cm} $\rank \rn( \A^{(n)}) \ \ =\ \  1\,$.
}
\vspace{4mm}

\noindent
{\bf Remark 1.3} ~The notion of height of a Hecke
algebra representation
was introduced in \cite{Gur} for the special case of
the representations generated by constant $R$-matrices.
There it was named {\em the rank of the $R$-matrix}.
We have changed the
name
here in order to avoid a possible
confusion with the standard notion of $\rank$ of a matrix.
Note that the use of
the
term `$\height$' is suggested by
the fact that imposing condition (\ref{a8}) for the representation $\rx$
results in vanishing of any central (and primitive) idempotent related to
a Young diagram (standard tableaux) containing
more then $n$ boxes in one of its columns.
\vspace{2mm}

\noindent
{\bf Remark 1.4} ~In view of Remark 1.1,
the whole sets of $q$-antisymmetrizers \\
$\{ \rx(\A^{(i,n+i)}) \}_{i=1,\dots ,k-n}$~
and $\{ \rx(\A^{(j,n+j-1)}) \}_{j=1,\dots ,k-n+1}$~
satisfy conditions (\ref{a8}) and (\ref{a9}), respectively.
\vspace{2mm}

\noindent
{\bf Remark 1.5} ~Instead of using (\ref{a8})
one can impose the condition
\be
\lb{11}
A^{(n+1)} = 0
\ee
at the algebraic level.
This is the way how generalized Temperley-Lieb-Martin algebras
are defined (cf. \cite{M}). Below we present several useful equivalent
forms of this condition.
\vspace{2mm}

\noindent
{\bf Lemma 1.1}
{\em
The condition (\ref{11})
is equivalent to any of the following relations
}
$$
\A^{(n)}\, g_n \dots g_2 g_1 =
(-1)^{n{-}1}q\,\q{n} \A^{(n)}\A^{(2,n+1)}\ , \eqno\mbox{(1.12a)}
$$
$$
g_1 g_2 \dots g_n\A^{(n)} =
(-1)^{n{-}1}q\,\q{n} \A^{(2,n+1)}\A^{(n)}\ , \eqno\mbox{(1.12b)}
$$
$$
g_n \dots g_2 g_1\A^{(2,n+1)} =
(-1)^{n{-}1}q\,\q{n} \A^{(n)}\A^{(2,n+1)}\ , \eqno\mbox{(1.12c)}
$$
$$
\A^{(2,n+1)} g_1 g_2 \dots g_n =
(-1)^{n{-}1}q\,\q{n} \A^{(2,n+1)}\A^{(n)}\ , \eqno\mbox{(1.12d)}
$$
$$
\A^{(n)}\A^{(2,n+1)}\A^{(n)} =
\q{n}^{-2} \A^{(n)}\ , \eqno\mbox{(1.12e)}
$$
$$
\A^{(2,n+1)}\A^{(n)}\A^{(2,n+1)} =
\q{n}^{-2} \A^{(2,n+1)}\ . \eqno\mbox{(1.12f)}
$$

\setcounter{equation}{12}
\renewcommand\theequation{\thesection.\arabic{equation}}

\noindent
{\bf Proof.}
Applying repeatedly (\ref{a4}) for the $q$-antisymmetrizers that appear as
last factors in the resulting products and using
(\ref{a6}), (\ref{a7}) we find
\ba
\nn
\A^{(n{+}1)} &=&
{1\over \q{n+1}} \A^{(n)} (q^n-\q{n}g_n)\A^{(n)}
\\
\nn
&=& {1\over \q{n+1}}\left\{	 q^n \A^{(n)} - \A^{(n)} g_n
(q^{n{-}1}	-\q{n-1}g_{n{-}1})\A^{(n{-}1)}\right\}
\\
\nn
 &=& {1\over \q{n+1}}\left\{\A^{(n)}(q^n - q^{n-1}g_n) +
\A^{(n)}g_n g_{n{-}1}(q^{n{-}2}-\q{n-2}g_{n-2})\A^{(n-2)}\right\}
\\
\lb{alt}
\dots &=&
{1\over \q{n+1}}\A^{(n)}\left(q^n-q^{n{-}1}g_n+
\dots + (-1)^n g_n g_{n-1}\dots g_1\right)\ .
\ea

Next, we apply $\A^{(2,n{+}1)}$ to the both sides of Eq.(\ref{alt}).
Using again (\ref{a6}) and (\ref{a7}) we obtain
$$
\A^{(n{+}1)} =
{1\over \q{n+1}}\A^{(n)}\left\{
q\q{n} \A^{(2,n+1)}+(-1)^n g_n g_{n-1}\dots g_1\A^{(2,n+1)}\right\}\ .
$$
Taking into account the relation
$
(g_n g_{n{-}1}\dots g_1) \A^{(2,n+1)}= \A^{(n)}(g_n g_{n{-}1}\dots g_1)
$
which is a consequence of (\ref{a1}) we end up with
\be
\lb{alt2}
\A^{(n{+}1)} =
{1\over \q{n+1}}\left\{ q\q{n}\A^{(n)} \A^{(2,n+1)}+
(-1)^n \A^{(n)}g_n g_{n-1}\dots g_1\right\}\ .
\ee
This proves the equivalence of Eqs.(\ref{11}) and (1.12a).
A similar argument using iteratively a substitution of
the first $q$-antisymmetrizer in the right hand side of (\ref{a4}) implies the
equivalence of (\ref{11}) and (1.12b). Condition (\ref{11}) is
transformed to
the forms (1.12c) and (1.12d) in the same manner starting from Eq.(\ref{a5}).

To show equivalence of (\ref{11}) to (1.12e) and to (1.12f)
one should employ Eqs.(\ref{a4}) and(\ref{a5}), respectively.
We shall treat the case of (1.12e) here.

Consider the difference
\ba
\nn
&&
\q{n}^2\A^{(n)}\A^{(2,n+1)}\A^{(n)} - \A^{(n)}
\, =\, \A^{(n)}\left(\q{n}^2\A^{(2,n+1)}-1\right)\A^{(n)}
\\
\nn
&=&
\A^{(n)}\left\{\q{n}\A^{(2,n)}(q^{n-1}-\q{n-1}g_n)\A^{(2,n)}
-1\right\}\A^{(n)}
\\
\nn
&=&\q{n-1}\A^{(n)}(q^n-g_n)\A^{(n)}
= \q{n-1}\q{n+1}\A^{(n+1)}\ ,
\ea
where we have again used the definition (\ref{a4}) and the relations
(\ref{a7}).
Comparing the first and the last
lines of the calculation we deduce the equivalence
of conditions (\ref{11}) and (1.12e).
\eod
\vspace{2mm}

Eqs.(1.12a--f) display properties of the $\rank 1$
idempotents $\rx (\A^{(n)})$ that are hidden in (\ref{11}).
In fact, they are the basic technical tools
which one needs to effectively deal
with the $\height n$ Hecke algebra representations.

In the rest of the paper we make use of a special type of
representations of the algebras ${\cal H}_k(q)$
for which the representation space is given by $k$-th tensor
power of an ($n$-dimensional, in the case of $SL(n)$)
vector space $V$: $W=V^{\otimes k}$. These representations
are generated by constant or dynamical $R$-matrices of Hecke type.
The representations we are dealing with have the specific feature
that their $\height$, when defined
(i.e., for $k\ge n\,$), coincides with the dimension of the space
$V~$\footnote{In general, this need not be the case.
Examples with $\height\rho\neq\dim V$ were constructed in \cite{Gur}.
}.
Below we first illustrate the general notions introduced above on
the well known case of constant $SL(n)$-type $R$-matrices
relegating the study of dynamical $R$-matrices to Sections 3 and 4.

\vspace{0.5 cm}

\section{Representations generated by a constant $R$-matrix of $SL(n)$ type}
\setcounter{equation}{0}
\renewcommand\theequation{\thesection.\arabic{equation}}

\medskip

The $R$-matrix corresponding
to the Drinfeld-Jimbo deformation of
$SL(n)$ \cite{D,J} is an operator acting in
a tensor square of an $n$-dimensional vector space $V$
and given by
\be
\lb{Jimbo}
\R^{\a_1\a_2}_{\b_1\b_2} =
q^{\d_{\a_1\a_2}}\, \d^{\a_1}_{\b_2} \d^{\a_2}_{\b_1}\, + \,
(q - \bq)\,\theta_{\a_2\a_1}\, \d^{\a_1}_{\b_1} \d^{\a_2}_{\b_2}\ ,
\ee
(no summation in the right hand side) where the indices $\a$, $\b$ take
values from $1$ to $n$, and
$\theta_{\a\b} = \{ 1\; {\rm if} \; \a >\b\ , \;
0 \; {\rm if} \; \a \leq \b\ \}$.

This $R$-matrix is a particular representative of a family
of constant Hecke $R$-mat\-ri\-ces,
-
i.e., it
satisfies the braid relation and the Hecke
condition
\ba
\lb{yb}
&&\R_{12}\R_{23}\R_{12} = \R_{23}\R_{12}\R_{23}\ ,
\\
\lb{hecke}
&&\R^2 = \id + (q - \bq)\, \R\ .
\ea

Eqs.(\ref{yb}),(\ref{hecke})
imply that the matrices $\R_{12}$, $\R_{23}$ generate a representation of
${\cal H}_3(q)$ in $V^{\otimes 3}$. For an arbitrary $k$ the representation
$\rho_{\R,k} :\, {\cal H}_k(q)\to \aut(V^{\otimes k})$
generated by a constant Hecke $R$-matrix is defined by
\be
\lb{const}
\rho_{\R,k}(g_i)	= \R_{i \, i{+}1}
\ .
\ee

For representations $\rho_{\R,k}$
generated by the $R$-matrix (\ref{Jimbo}) we have
$$\height\rho_{\R,k} \,= \, n\quad {\rm if} \quad k\geq n\,.$$
The $\rank 1$
$q$-antisymmetrizers $\rho_{\R,k}(\A^{(i,n+i-1)})$ are most conveniently described
in terms of $q$-analogues of (co- and contravariant) Levi--Civit\`a tensors
which are solutions of the equations
\ba
\nn
\R_{\scriptscriptstyle \b_i \b_{i{+}1}}^{\scriptscriptstyle \a_i \a_{i{+}1}}
\e^{\scriptscriptstyle \a_1 \dots \b_i \b_{i{+}1}\dots \a_n} =
-{\bq}\, \e^{\scriptscriptstyle \a_1\dots \a_i \a_{i{+}1}\dots \a_n}\, , &&
\\
\nn
\e_{\scriptscriptstyle \a_1 \dots \b_i \b_{i{+}1}\dots \a_n}
\R_{\scriptscriptstyle \a_i \a_{i{+}1}}^{\scriptscriptstyle \b_i \b_{i{+}1}}
= -{\bq}\, \e_{\scriptscriptstyle \a_1\dots \a_i \a_{i{+}1}\dots \a_n}\, ,
&\quad i=1,2,\dots ,n{-}1\, .&
\ea
It is straightforward to prove that these equations have unique (up to
normalization) solutions. The $\rank 1$ condition (\ref{a9}) follows as a
corollary.

In the special case of representations $\rho_{\R,k}$
generated by the $R$-matrix (\ref{Jimbo})
the only nonvanishing components of the $\e$-tensors have
pairwise different indices $\a_1,\a_2,\dots \a_n\,$, and can be chosen as
\be
\lb{e}
\e^{\scriptscriptstyle \a_1 \a_2 \dots \a_n} =
\bq^{n(n-1)/2}\, (-q)^{\ell(\sigma)}\ , \qquad
\e_{\scriptscriptstyle \a_1 \a_2 \dots \a_n} =
(-q)^{\ell(\sigma)} \; .
\ee
Here $\ell(\sigma)$ is the length of the permutation
$\sigma = {1_{~}, 2_{~}, \dots ,n_{~}\choose \a_1, \a_2, \dots ,\a_n}$.

\vspace{2mm}

The $\rank 1$ $q$-antisymmetrizers are expressed in terms of the $\e$-tensors
as
\be
\lb{A}
\rho_{\R,k}(A^{(i,n+i-1)}) =
{1\over [n]!}\, \eupi{i}\edoi{i}\ , \quad i=1,\dots ,k-n+1\ .
\ee
Here (by analogy with the matrix notation)	we substitute
the vector space indices of $\e$-tensors by their labels:
$\a_i \to i$.
The "bra" and "ket" notation of $\e$-tensor
indices is used in order to distinguish
labels of matrix spaces from those of vector spaces.
One should have in mind the following symbolic decomposition for
the matrix space label:
$i = \raisebox{2pt}{$ |i{\cal i}$}\otimes\raisebox{-2pt}{${\cal h}i|$}$.
For example, the equation
$A_i\, u^{|i{\cal i}}
\left(\equiv A^{|i{\cal i}}_{\;\;\; {\cal h}i|}\, u^{|i{\cal i}} \right)
= v^{|i{\cal i}}$
is to be understood as
$\sum_{\b_i} \, A^{\a_i}_{\;\b_i} \, u^{\b_i} = v^{\a_i}$.

Finally, we shall adapt for $\rho_{\R,k}$
those formulas (1.12a-f) which will be used
in Section 5.
Written in terms of the $\e$-tensors
the relations (1.12b), (1.12c), (1.12e) and (1.12f) assume the form
\be
\lb{R11b}
\rho_{\R,n{+}1}(g_1\dots g_n)\eup \equiv
\R_{12}\dots \R_{n \, n{+}1}\eup = q\, \eu2\, \N{1}{n{+}1}\ ,
\ee
\be
\lb{R11c}
\rho_{\R,n{+}1}(g_n\dots g_1)\eu2 \equiv
\R_{n \, n{+}1}\dots\R_{12}\eu2 = q\,\eup\, \K{n{+}1}{1}\ ,
\ee
\be
\lb{KN}
K\, N =  N\, K	= \id\ .
\ee
Here the matrices $N$ and $K$ are defined as
\ba
\lb{N}
\N{1}{n{+}1} &=& {(-1)^{n-1}\over \q{n-1}!}\, \ed2 \eup\ ,
\\
\lb{K}
\K{n{+}1}{1} &=& {(-1)^{n-1}\over \q{n-1}!}\,	\edo \eu2\ ,
\ea
and for the $\e$-tensors given by Eq.(\ref{e}) we have just
\be
\lb{nkid}
N = K = \id\ .
\ee

\vspace{0.5 cm}

\section{$SL(n)$-type  dynamical $R$-matrices}

\setcounter{equation}{0}
\renewcommand\theequation{\thesection.1\alph{equation}}

\medskip

We now turn to the dynamical $R$-matrix defined in the Introduction.
In order to present the QDYBE in a form suitable for our purposes we
shall introduce a set of commutative variables
\be
\lb{defX}
X^i\;,\,i=1,\ldots , n,\quad [X^i ,\, X^j ] = 0\,,\quad
\prod\limits_{i=1}^n X^i = 1
\ee
which play the role of elementary shift operators for $p_i\,:$
\be
\lb{0.8}
p\, X^i\, =\, X^i\, (\, p \,+\, v^{(i)}\, )\,.
\ee
The elements  $X^i\;$ and $q^{p_i}\;$ provide a realization of
(the Weyl's form of) the
canonical commutation relations. Note that in concrete
applications they can be naturally
identified with (a subset of) dynamical variables of a model
(see, e.g., \cite{AF}).

\setcounter{equation}{1}
\renewcommand\theequation{\thesection.\arabic{equation}}

Let us arrange the auxiliary variables $X^{i}$
into a unimodular diagonal matrix,
\be
\lb{X}
X = diag\{X^1, \dots ,X^n\}\ , \qquad \det (X) = 1\ .
\ee

The {\em Hecke-type dynamical $R$-matrix} is characterized by the
following set of relations:
\be
\lb{dybe}
\DR{12} (X_1 \DR{23} X_1^{-1}) \DR{12} =
(X_1 \DR{23} X_1^{-1}) \DR{12} (X_1 \DR{23} X_1^{-1})\,,
\ee
\ba
\lb{dhecke}
\Rp^2 &=& \id + (q- \bq)\Rp\,,
\\
\lb{add}
\DR{12} X_1 X_2 &=& X_1 X_2 \DR{12}\,.
\ea
Here the first and second relations are
the dynamical Yang-Baxter equation and the Hecke condition, respectively.
A condition of type (\ref{add}),
although not always imposed on dynamical $R$-matrices,
is also necessary in our treatment. As we shall see below,
it ensures that conditions (\ref{a3}) for the Hecke
algebra representations generated by $\Rp$ are satisfied.
\medskip

Following \cite{I2} we shall consider dynamical $R$-matrices of the form
\be
\lb{Rp}
\DDR{i_1 i_2}{j_1 j_2} =
a_{i_1 i_2}(p) \, \d^{i_1}_{j_2} \d^{i_2}_{j_1}\, + \,
b_{i_1 i_2}(p)\, \d^{i_1}_{j_1} \d^{i_2}_{j_2}\ ,
\qquad i_{1,2} \, , \, j_{1,2} = 1,\dots ,n
\ee
(there is no summation over repeated upper and lower indices
in the right hand side); in order to have a unique decomposition in terms of the
unit and the permutation matrices in the tensor square of spaces
we impose the condition $b_{ii}(p)=0\,$.
This special class of $p$-dependent Hecke $R$-matrices will be called
{\em dynamical $R$-matrices of an $SL(n)$-type~}\footnote{
In \cite{EV} these $R$-matrices were called  $GL(n)$-type
$R$-matrices. We call them $SL(n)$-type $R$-matrices instead
since we impose the additional condition $\sum_{i=1}^n p_i = 0$
(see (\ref{0.2})) on the variables $p_i$.
}.

The unknown functions $a_{ij}(p)$, $b_{ij}(p)$
in the  Ansatz (\ref{Rp}) are to be fixed by the conditions
(\ref{dybe}-\ref{add}). The Hecke condition
(\ref{dhecke}) gives
\be
\lb{btt}
b_{ii} = 0\ ,\ b_{ij}+ b_{ji}= q - \bq \; ,
\quad \mbox{\rm for}\quad i\neq j \ ,
\ee
\be
\lb{att}
a_{ij} \, a_{ji} - b_{ij} \, b_{ji} = 1 \; ,
\quad \mbox{\rm for}\quad i\neq j \ ,
\ee
\be
\lb{aii}
a_{ii}^{2} - (q-\bar q) \, a_{ii} = 1 \; .
\ee
The last equation has two solutions: $a_{ii} = \pm\, q^{\pm 1}$
for each $i$.
Below we consider only the case $a_{ii} =q,\; \forall i$ (the other
cases correspond, in particular, to quantum supergroups
and have been considered in \cite{I2}).
Finally, the dynamical Yang-Baxter equation (\ref{dybe}) and Eq.(\ref{add})
impose the constraints
$$
a_{ij}(p_1, \dots , p_n) = a_{ij}(p_{ij}) \;\; , \;\;\;
b_{ij}(p_1, \dots , p_n) = b_{ij}(p_{ij}) \; ,
$$
\be
\lb{btt2}
 b_{ij} \, b_{jk} \, b_{ki} + b_{ik} \, b_{kj} \, b_{ji} = 0 \; ,
\ee
\be
\lb{fff}
b_{ij}(p_{ij} + 1) =
\frac{b_{ij}(p_{ij}) \, q}{\bar q + b_{ij}(p_{ij})} \; ,
\ee
where $p_{ij} := p_i - p_j$.
For $a_{ii} =q$
the general solution of (\ref{btt}) - (\ref{fff})
can be written as
\cite{I2}
\be
\lb{ab}
a_{ij}(p) = \a_{ij}(p_{ij}) \xi_{ij}(p_{ij})\ , \quad
b_{ij}(p) = q - \xi_{ij}(p_{ij})\ ,
\ee
where
$\xi_{ij}(p)$	are expressed as the following ratios:
\be
\lb{xi}
\xi_{ij}(p)	= {f(p_{ij}-1,\b_{ij})\over f(p_{ij},\b_{ij})}\ , \quad
f(p,\b)=\bq^p + \q{p}\b\ .
\ee
Here $\beta_{ij}(p_{ij}) =\beta_{ij}(p_{ij}+1)\,$.
We shall consider $\beta_{ij}$ as constant parameters
since their functional dependence does not change any of the results below.
The function $f(p,\b)$
satisfies the finite difference equation
\be
\lb{f}
f(p+1,\b) + f(p-1,\b) = \q{2}f(p,\b)
\ee
with the initial conditions
\be
\lb{initial}
f(0,\b) = 1\ , \qquad f(1,\b)=\bq + \b\ .
\ee
Equations (\ref{f}) can be deduced from
(\ref{fff}).

Not all of the remaining in (\ref{ab}), (\ref{xi})
parameters $\a_{ij}(p_{ij})$ and $\b_{ij}$
are independent. The relations between them are given by
\be
\lb{al}
\a_{ii} = 1\,,\quad \a_{ij}(p_{ij})\a_{ji}(p_{ji}) = 1\ ,
\ee
\setcounter{equation}{0}
\renewcommand\theequation{\thesection.17\alph{equation}}
\be
\lb{bt}
\b_{ii} = 0\,,\quad  \b_{ij}+\b_{ji}= q - \bq\
\quad \mbox{\rm for}\quad i\neq j \ ,
\ee
\be
\lb{bt2}
\b_{ij}\b_{jk}\b_{ki}+\b_{ik}\b_{kj}\b_{ji} = 0 \; .
\ee

\setcounter{equation}{0}
\renewcommand\theequation{\thesection.18\alph{equation}}

An easy way to solve Eqs.(\ref{bt}), (\ref{bt2})
is to make the substitution
\be
\lb{pi1}
\b_{ij} = \frac{q-\bq}{1 - \pi_{ij}} \quad \Leftrightarrow \quad
\pi_{ij} = \frac{\b_{ij} - q + \bq}{\b_{ij}} \; , \quad
\mbox{\rm for}\ i\neq j\ .
\ee
We stress that the parameters $\pi_{ii}$ are not fixed here
and can be chosen as arbitrary constants.
In terms of
the
new variables $\pi_{ij}$ equations (\ref{bt}), (\ref{bt2})
take the simple form
\be
\lb{pi2}
\pi_{ij} \, \pi_{ji} = 1 \; , \quad
\pi_{ij} \, \pi_{jk} \, \pi_{ki} = 1 \; ,
\ee
and are solved by
$\pi_{ij} = \prod_{k=i}^{j{-}1}\pi_{k \, k{+}1} =\pi_{ji}^{-1}$, for $i<j$.
Hence,

\setcounter{equation}{18}
\renewcommand\theequation{\thesection.\arabic{equation}}

\be
\lb{betas}
\b_{ij}	=
{(q- \bq) \prod_{k=i}^{j-1} \, \b_k \over
\prod_{k=i}^{j-1} \, \b_k - \prod_{k=i}^{j-1} \,
(\b_k - q +\bq)}\ , \quad
\mbox{\rm for} \quad i<j\ ,
\ee
where the remaining $(n-1)$ parameters
$\b_i \equiv \b_{i \, i+1}$ are independent.
\vspace{2mm}

We shall describe two types of transformations on the set
of $SL(n)$-type dynamical $R$-matrices.
\vspace{2mm}

\begin{enumerate}
\item
a simple version of {\em twisting} for dynamical $R$-matrices:
$$
{\Rp}_{21} \to \xR{F}_{21} = \hF \, \Rp \, \hF^{-1} \; ,
$$
where
$\hF(p) = F_{12}(p) P_{12}$ and
\be
\lb{twi}
F_{12} \equiv F^{i_1 i_2}_{j_1 j_2}(p)	=
\d^{i_1}_{j_1} \, \d^{i_2}_{j_2} \, \psi_{i_1 i_2}(p)  \; ,
\ee
(an analog of Drinfeld's twist, see \cite{Drinfel}).
An explanations on how this twist works is given in
the
two Lemmas below.
\vspace{2mm}

\noindent
{\bf Lemma 3.1} {\em Let $\hF(p)$ be an operator acting
in $V \otimes V$. If
\be
\hF^{-1}_{12} \, X_1^{-1} \, \hF_{23} =
\hA_{123} \, X_3^{-1} \, \hA_{123}
\lb{tf}
\ee
for some operator $\hA(p)$ acting in $V \otimes V \otimes V$ and
\be
\Rp_{12} \, \hA_{123}=\hA_{123} \, \Rp_{23}
\lb{ta}
\ee
then the matrix $\xR{F}_{21}=\hF \Rp \hF^{-1}$ satisfies
(a version of -- see the proof)
the QDYBE (\ref{dybe}).}
\vspace{2mm}

\noindent
{\bf Proof.} Substitute $\xR{F}_{21}$ in the QDYBE taking into
account that the QDYBE has two
more equivalent forms:
$$
\Rp_{23} \, X_3^{-1} \, \Rp_{12} \, X_3 \, \Rp_{23} =
X_3^{-1} \, \Rp_{12} \, X_3 \, \Rp_{23} \,
X_3^{-1} \, \Rp_{12} \, X_3 \; ,
$$
$$
\Rp_{21} \, X^{-1}_1 \Rp_{32} X_{1} \, \Rp_{21} =
X^{-1}_1 \Rp_{32} X_{1} \, \Rp_{21} \, X^{-1}_1 \Rp_{32} X_{1} \; .
$$
The first
equation
results from repeated application of (\ref{add}) to
(\ref{dybe}). The second
equation is obtained from the first one by simply
permuting the subscripts $1$ and $3$.
\eod
\vspace{2mm}

\noindent
{\bf Lemma 3.2} {\em
Let $F_{12}$ be
the
diagonal matrix (\ref{twi}) where
$\psi_{ij}\, \psi_{ji}=1$ and $\psi_{ii}=1$.
Assume that $\psi_{ij}$ depends on $p_{ij}$ only.
Then,
Eqs.(\ref{tf}),
(\ref{ta}) are satisfied by
$\hA =A P_{23}P_{12}$ where the matrix $A$ is
diagonal, $A^{ijk}_{abc}=a_{ijk}\d^i_a\d^j_b\d^k_c \,,$
the elements $a_{ijk}$ being given by }
\be
a_{ijk}=\left\{
\begin{array}{ll}
\, \psi_{ik} \, \psi_{jk} & {\rm if}\ i\neq j\ ,\\
\, [\psi_{ik}(p_{ik}+1)]^2 & {\rm if}\ i=j\ .
\end{array}
\right.
\lb{sa}
\ee

\noindent
{\bf Proof.} The operator $A_{123}$ is symmetric in the first
two indices, $A_{123}=A_{213}$ which implies that it commutes with
any $R$-matrix of the form (\ref{Rp}).
Therefore, (\ref{ta}) is satisfied.

Eq.(\ref{tf})
can be checked directly.
\eod

These Lemmas demonstrate that
the
operator
(\ref{twi})
indeed generates a twist leading to the changes
$\a_{ij} \rightarrow \a_{ij} \psi^{2}_{ji}$, $\b_{ij} \to \b_{ij}$
of the parameters in (\ref{ab}), (\ref{xi}).

\item
{\em canonical transformations}
of the dynamical parameters $p_i \to p_i + c_i$,
where $c_i,\; i=1,\dots ,n$ are arbitrary constants satisfying
the
condition $\sum_{i=1}^n c_i = 0$.
\end{enumerate}

\vspace{2mm}
We conclude the Section by
a brief discussion of the structure of
the family of $SL(n)$-type
dynamical $R$-matrices (\ref{Rp}), (\ref{ab}), (\ref{xi}).
There are two essentially different domains for
the parameters $\b_i$ of this family.
\vspace{2mm}

{\bf a)}
$\b_i \neq 0$ and $\b_i \neq q - \bq$, for all $i$.

In this case the whole family
(\ref{Rp}), (\ref{ab}), (\ref{xi}) can be generated starting from
any particular representative with the use of
the two types of transformations described above.

Indeed, the parameters $\a_{ij}$ can be excluded with the help of
a twist. Then,
performing a canonical transformation of the form
$$
q^{2p_{ij}} \rightarrow
q^{2p_{ij}} \, \pi_{ij} =
q^{2p_{ij}} \, \prod_{k=i}^{j-1} \,
\frac{\b_k - q + \bq}{\b_k} \;\; {\rm for} \;\; i<j \; ,
$$
for instance, one
excludes the parameters $\b_{ij}$ from the
Ansatz (\ref{Rp}), (\ref{ab}), (\ref{xi})
and passes to
a
dynamical $R$-matrix with
\be
\lb{xi2}
\xi_{ij}(p) = {\q{p_{ij}-1} \over \q{p_{ij}}}
\ee
(cf.
(\ref{xi}) -- (\ref{bt2})).
This $R$-matrix
is the limiting case $\b_i \to \infty$ of our family,
and it is this type
of
dynamical $R$-matrix which is discussed in
\cite{AF}.

{\bf b)}
Either all $\b_{i} = q-\bq$, or all $\b_i=0$.\\
\noindent
We shall consider the first case $\b_i=q-\bq$ for which
\be
\lb{case1}
\b_{ij} = \left\{
\begin{array}{lll}
q-\bq &\mbox{\rm for}& i<j
\\
0 & \mbox{\rm for}& i\geq j
\end{array}
\right.
\ee
and put
$\a_{ij}(p_{ij})= \mbox{\rm const}_{ij}$.
In this case the $R$-matrix (\ref{Rp}), (\ref{ab}), (\ref{xi})
becomes independent
on
the dynamical variables $p_i$ and is reduced
to the constant $R$-matrix
describing the multiparametric \cite{Man} deformations of $GL(n)$
which are all twist-equivalent.

With the particular choice
\be
\lb{case}
\a_{ij} = \left\{
\begin{array}{ccc}
q&\mbox{\rm for}&i<j
\\
1&\mbox{\rm for}&i=j
\\
\bq&\mbox{\rm for}&i>j
\end{array}
\right.
\ee
one reproduces the standard $SL(n)$-type $R$-matrix
(\ref{Jimbo}).
\vspace{2mm}

\noindent
{\bf Remark 3.1}
In the intermediate cases (where only a
part of the parameters $\beta_{i}$
are equal to $0$ or $q-\bq$) the corresponding dynamical $R$-matrix $\Rp$
contains the (dynamical and constant) $R$-matrices described in {\bf a)}
and {\bf b)}
as submatrices.
\vspace{2mm}

\section{
Representations generated by $SL(n)$-type  dynamical $R$-matrices}

\setcounter{equation}{0}
\renewcommand\theequation{\thesection.1\alph{equation}}
\medskip

Now we are in a position to
introduce the Hecke algebra representations
associated with Hecke-type dynamical $R$-matrices.
\vspace{3mm}

{\samepage
\noindent
{\bf Proposition 4.1} {\em
Let $\Rp$ be a dynamical $R$-matrix of the Hecke type.
The matrices
\be
\lb{drep1}
\rho_{\Rp,k}(g_{i}) =
(X_1 X_2\dots X_{i{-}1})\DR{i,i{+}1}
(X_1 X_2\dots X_{i{-}1})^{-1} \; , i =1,\dots,k{-}1,
\ee
generate a Hecke algebra representation,
$\rho_{\Rp,k}: {\cal H}_k(q) \to \aut(V^{\otimes k})$.}
}
\vspace{2mm}

\noindent
{\bf Proof.} Obviously, equation
(\ref{dybe}) imply that the matrices
$\rho_{\Rp,k}(g_{i})$ and $\rho_{\Rp,k}(g_{i+1})$
satisfy the braid relations (\ref{a1}). Then, the
conditions (\ref{add}) ensure that
the matrices (\ref{drep1}) satisfy (\ref{a3}) and, therefore,
(\ref{drep1}) represent the generators
of the braid algebra ${\cal B}_k$. Finally,
the Hecke conditions (\ref{a2})	for the generators
(\ref{drep1}) follow from the Hecke property
(\ref{dhecke}) of the dynamical $R$-matrix
(\ref{a2}).
\eod
\vspace{3mm}

\noindent
{\bf Remark 4.1}
In contrast with the case of
constant Hecke $R$-matrix (\ref{const})
the representations generated by a dynamical Hecke $R$-matrix
are nonlocal; in other words, the matrices $\rho_{\Rp,k}(g_i)$
act nontrivially as diagonal matrices
on $V_j \,$ with $j<i\,$ (and not merely on $V_i\otimes V_{i+1}$). Only the
representation of the first generator with $i=1$ has the
usual `locality' property.

\vspace{2mm}
\noindent
{\bf Remark 4.2} One can construct representations
equivalent to $\rho_{\Rp,k}$
in which some other generator is `localized' instead $g_1\;$. For
instance, the representation $\overline{\rho}_{\Rp,k}$	which localizes
$\overline{\rho}_{\Rp,k}(g_{k{-}1})$ is given by
\be
\lb{drep1-2}
\overline{\rho}_{\Rp,k}(t) =
(X_1 X_2\dots X_k)^{-1}\rho_{\Rp,k}(t)(X_1 X_2\dots X_k)\ ,	\quad
\forall	t\in {\cal H}_k(q)
\ee
so that
\be
\lb{drep2}
\overline{\rho}_{\Rp,k}(g_i) =
(X_{i{+}2}\dots X_k)^{-1}\DR{i,i{+}1} (X_{i{+}2}\dots X_k)
\ .
\ee
Note that in addition to the nonlocal property
the representation matrices of $\overline{\rho}_{\Rp,k}$
depend explicitly on $k$.
\vspace{2mm}

\setcounter{equation}{1}
\renewcommand\theequation{\thesection.\arabic{equation}}

From now on we shall restrict ourselves to discussing
those representations
$\rho_{\Rp,k}$ which are generated by the $SL(n)$-type dynamical
$R$-matrices (\ref{Rp}), (\ref{ab}), (\ref{xi}).
For $k\geq n$ all these representations
are of $\height n$.
The $\rank 1$ $q$-antisymmetrizers are conveniently expressed in terms of
dynamical $\E$-tensors $\Eup$ and $\Edo$, which are the unique
(up to normalization)
solutions of the equations
\be
\lb{E}
\begin{array}{rcl}
\rho_{\Rp,k}(g_i)\, \Eup &=& - \bq\: \Eup\ ,
\\
\Edo \,
\rho_{\Rp,k}(g_i) &=& - \bq\:  \Edo\ ,
\end{array}
\quad  1\leq i\leq n-1\ .
\ee
The only nonvanishing components of these $\E$-tensors
have pairwise different indices $i_1,i_2,\dots ,i_n$ and look like
\ba
\lb{Eup}
\E^{i_1 i_2 \dots i_n}(p) &=& (-1)^{\ell(\sigma)}\!\!
\prod_{(j,i)\in J(\sigma)}\!\!\!\! \a_{ji}(p_{ji})\,
\prod_{1\leq a<b\leq n}	\!\!\!\!\xi_{i_a i_b}(p_{i_a i_b})\ ,
\\
\lb{Edo}
\E_{i_1 i_2 \dots i_n}(p) &=& (-1)^{\ell(\sigma)} \!\!
\prod_{(j,i)\in J(\sigma)} \!\!\!\!\a_{ij}(p_{ij})\ .
\ea
Here $\ell(\sigma)$ is the length of the permutation
$\sigma = {1_{~}, 2_{~},\dots , n_{~}\choose i_1, i_2, \dots ,i_n}$,
and
$$
J(\sigma) = \{ (i_a,i_b): a<b, i_a>i_b \}\ .
$$
The dynamical $\E$-tensors (\ref{Eup}), (\ref{Edo})
are normalized so that they would
coincide with the constant $\e$-tensors (\ref{e})
in the case (\ref{case1}), (\ref{case}).

Now the expressions for $\rank 1$ $q$-antisymmetrizers
in the representations $\rho_{\Rp,k}$
are given by
\ba
\nn
\lefteqn{
\rho_{\Rp,k}(\A^{(i,n+i-1)}) =
}
\\
\lb{DA1}
&&
{1\over \q{n}!}\,
( X_1\dots X_{i{-}1})
\Eupi{i}\Edoi{i}
( X_1\dots X_{i{-}1})^{-1}
\ .
\ea
The
numerical
coefficient in this formula is calculated with the
use of the relation
\be
\lb{norm}
\Edo \Eup =  \q{n}!\ ,
\ee
which is proved in the Appendix.

We conclude the discussion on dynamical $R$-matrices
by writing down formulas (1.12a), (1.12d),
and (1.12e), (1.12f) for the representation
$\rho_{\Rp,n+1}\;$:
\be
\lb{DeR}
\Edo\, \rho_{\Rp,n+1}(g_n g_{n{-}1}\dots g_1) =
q\, \K{n{+}1}{1}(p)\,	X_1\Ed2 X_1^{-1}\ ,
\ee
\be
\lb{DeR2}
X_1 \Ed2 X_1^{-1} \rho_{\Rp,n+1}(g_1 g_2\dots g_n) = q\,
\N{1}{n{+}1}(p)\, \Edo\ ,
\ee
\be
\lb{DKN}
K(p)\, N(p) =  N(p)\, K(p)	=
\id\ .
\ee
Here the matrices $N(p)$, $K(p)$
are defined as
\ba
\lb{Np}
\N{1}{n{+}1}(p) &=&
{(-1)^{n-1}\over \q{n-1}!}\, X_1\Ed2 X_1^{-1} \Eup\ ,
\\
\lb{Kp}
\K{n{+}1}{1}(p) &=&
{(-1)^{n-1}\over \q{n-1}!}\,	\Edo X_1 \Eu2 X_1^{-1}
\ .
\ea
For the $SL(n)$-type dynamical $R$-matrices
the matrices $N(p)$, $K(p)$
are diagonal.
Inserting formulas (\ref{Eup}), (\ref{Edo}) for the dynamical
$\e$-tensors into (\ref{Np}), (\ref{Kp}) and using (\ref{norm}),
one ends up with the following expressions for their diagonal components:
\ba
\lb{npi}
N^{\, i}_{\: i}(p)
\, = \,  (K^{\, i}_{\: i}(p))^{-1}
&=& \prod_{j\neq i} \, \a_{ij} \, (p_{ij} - \theta_{ji}) \,
\xi_{ij}(p_{ij})
\ ,
\ea
where
$\theta_{ji} = \{ 1\; {\rm if} \; j>i\ , \;
0 \; {\rm if} \; j<i\ \}$.

\vspace{0.5 cm}

\section{
Quantum matrix algebra ${\cal A}(\Rp, \R)$:
quantum determinant and inversion formula}
\setcounter{equation}{0}
\renewcommand\theequation{\thesection.\arabic{equation}}

\medskip

We shall apply the above technique
to the quantum matrix algebra ${\cal A}$ which is defined as
follows (cf. Introduction).
\vspace{2mm}

{\samepage
\noindent
{\bf Definition 5.1~} {\em Let $\F$ be the field of the
complex meromorphic functions of the (commuting) variables
$p_j\,,j=1,\ldots ,n\,$.
Let $\Rp$ be a dynamical
$R$-matrix of an $SL(n)$-type and
$\R$ be a
constant $SL(m)$-type $R$-matrix.
Assume that both $\Rp$ and $\R$ satisfy
the Hecke condition (0.4) with the same value of $q$.
Then ${\cal A} = {\cal A}(\Rp,\R)\;$
is a complex algebra with $1$ that is generated	by $\F$ and the
$mn\;$ elements $a^i_{\a}$ ($i\; =
1,\ldots , n$ and $\a\; =1,\ldots ,m$),
satisfying the relations
\ba
\lb{aa}
&&\Rp_{12}\, \ai{1} \ai{2} \, = \,
\ai{1} \ai{2} \R_{12}\ ,
\\
\lb{ap}
&& a f(p) \;=\; X f(p) X^{-1} \; a\ , \qquad \forall f(p)\in\F\ ,
\ea
where $X\;$ is a unimodular diagonal matrix (\ref{X})
whose diagonal elements $X^i\;$ satisfy (3.1a), (3.1b).
}
}
\vspace{2mm}

\noindent
{\bf Remark 5.1~} The definition above is given for
arbitrary $m$ and $n$. However in the sequel we shall discuss
the case $m=n$ only.
\vspace{2mm}

\noindent
{\bf Remark 5.2~} For the applications envisaged
here, the field $\F$ of meromorphic functions of $p_j$
can be replaced by its
subfield of rational functions of $q^{p_j}$ (as stated in the
Introduction).
Then we should just require
$$
X^k q^{p_{ij}}(X^k)^{-1}= q^{p_{ij}+\delta_{jk}-\delta_{ik}}
$$
instead of (3.1b).
Note that for $q$ a root of unity, $q^{2h}=1$ (cf. Eq.(6.8) below)
$p_{ij}$
are only determined up to an additive integer multiple of $2h$.
\vspace{2mm}

\noindent
{\bf Remark 5.3~} More general matrix algebras are of interest in which
the $R$-matrices on both sides of the quadratic relations (\ref{aa}) are
allowed to depend on possibly different sets of commuting variables $p\;$
and $p'\;$
\be
\lb{QQ}
\Rp_{12}\, Q_{1} Q_{2} \, = \,
Q_{1} Q_{2} {\hat R}'(p')_{12}\ ,
\ee
while the shift properties assume the form
\ba
\lb{pQ}
 p_{kl} \, Q^i_j \, =\, Q^i_j\, ( p_{kl} + \d^i_k - \d^i_l ) &,&
{p'}_{kl} \, Q^i_j \, =\, Q^i_j\, ( {p'}_{kl} + \d_{jk} - \d_{jl} )\;.
\ea
Such ${\cal A}(\Rp,{\hat R}'(p') )\;$ can be treated in much the same way
or reduced to the study of two matrix algebras of the above type
setting $Q^i_j = a^i_{\a} {\overline a}^{\a}_j\;$,
where $a$ and $\overline a$ satisfy exchange relations of
the type (\ref{0.6}) and (\ref{0.12}), respectively (see \cite{FHIPT}).
Note that dynamical quantum groups (introduced in \cite{Fe}) are
defined by relations similar to (\ref{QQ}) and (\ref{pQ}) but with
the dynamical $R$- matrices (and momenta $p,p'$) related to each other
by some equivalence transformation
${\hat R}'(p') = X^{-1} \Rp X$.
Another desirable modification of the matrix algebra (\ref{aa})
corresponds
to the case when $\R$ is an $SO_q$ or $Sp_q$ constant
$R$-matrix.
In this case $\R\;$ and $\Rp\;$ satisfy a third order (Birman--Wenzl)
condition instead of the Hecke property (\ref{a2}) and the QDYBE (\ref{dybe})
have to be modified correspondingly.
\vspace{2mm}

\noindent
{\bf Remark 5.4~}
The algebra $\cal A$ differs from the one considered in \cite{CrGe} where the
counterpart of a matrix $b=X^{-1}a$
(denoted by $u^i_{\a}({\bar \omega})$ in Eq.(11) of
\cite{CrGe})
which commutes with $p$ is used for changing the basis of
chiral vertex operators. It is assumed in \cite{CrGe}
that the elements of $b$ only
depend on $p$
and hence commute
among themselves while in our case
this is not so. Indeed, the reflection equation subalgebra ${\cal M}(\R)$,
defined in Proposition 5.5 below, is non-commutative although its
elements commute with the $p$'s.
The difference is essential: as a result, Cremmer and
Gervais do not recover the standard (constant) $SL_q (n) \ R$-matrix for
$n>2$ but introduce instead new solutions of the Yang-Baxter
equation. One of the authors
(I.T.) would like to thank J.-L. Gervais and E. Cremmer for an enlightening
discussion on this point.
\vspace{2mm}

The term {\em "matrix algebra"~}\footnote{More conventional
quantum matrix algebras have been introduced in \cite{FRT,Man}
and their matrix nature were further investigated in \cite{Gur,GIOPS}.}
for the algebra ${\cal A} \equiv {\cal A}(\Rp, \R)\,$
is justified by the fact that we shall be able to (define and)
compute the determinant of $a\;$ -- as a function of $p\;$ -- and to find
the inverse of $a\;$. In the case of $2\times 2\;$ matrices the
determinant of $a\;$ was constructed in \cite{BF} (see also \cite{FeV})
for the special choice $\b_i\to\infty , \a_{ij}=1$ of the parameters.
We shall present the definition of the determinant
in a general setting.

\vspace{2mm}

\noindent
{\bf Definition 5.2} {\em Let $a=||a^i_{\a}||$ be the matrix
of generators of the algebra ${\cal A}(\Rp,\R)$. The determinant
of the matrix $a$ is given by
}
\be
\lb{deta}
\det (a)\, =\, {1\over \q{n}!}\,
\Edo \ai{1} \ai{2} \dots \ai{n} \eup\ .
\ee

The meaning of this definition is made clear by the following
three Propositions. The first and the third of them are
the quantum analogues of the basic determinant properties.
The second one allows to perform an $SL(n)$-reduction in
the algebra ${\cal A}(\Rp, \R)$.
\vspace{2mm}

\noindent
{\bf Proposition 5.1} {\em
The product $(\ai{1} \ai{2} \dots \ai{n})$ intertwines between
constant and dynamical $\e$-tensors:
}
\ba
\lb{ea}
\Edo \ai{1} \ai{2} \dots \ai{n} &=& \det (a) \, \edo \; ,
\\
\lb{ae}
\ai{1} \ai{2} \dots \ai{n} \eup &=& \Eup \det (a) \; .
\ea

\noindent
{\bf Proof.}
First, observe that due to the relations (\ref{aa}), (\ref{ap})
the product of $k$ matrices  $(\ai{1}\ai{2}\dots \ai{k})$
intertwines between the representations $\rho_{\R,k}$ and $\rho_{\Rp,k}$
of the algebra ${\cal H}_k(q)$.
Indeed,
\ba
\nn
&&(\ai{1}\dots \ai{k})\rho_{\R,k}(g_i) =\,
(\ai{1}\dots \ai{k})\R_{i(i{+}1)}
\\
\lb{intertwine}
= &&\ai{1}\dots\ai{i{-}1}
(\DR{i(i{+}1)}\ai{i}\ai{i{+}1})\ai{i{+}2}\dots \ai{k}
\\
\nn
= &&(X_1\dots X_{i{-}1})\DR{i(i{+}1)}(X_1\dots X_{i{-}1})^{-1}
(\ai{1}\dots \ai{k}) =\,
\rho_{\Rp,k}(g_i)(\ai{1}\dots \ai{k})\ .
\ea

In particular, one has
$$
(\ai{1}\dots\ai{n})\rho_{\R,n}(A^{(n)}) =
\rho_{\Rp,n}(A^{(n)})(\ai{1}\dots\ai{n})\ .
$$
Multiplying both sides by $\rho_{\R,n}(A^{(n)})$ from the right
or by $\rho_{\R,n}(A^{(n)})$ from the left and using
projector property of the $q$-antisymmetrizer one comes to
the equations
\ba
\lb{ea2}
\rho_{\Rp,n}(A^{(n)})(\ai{1}\dots\ai{n}) &=&
\rho_{\Rp,n}(A^{(n)})(\ai{1}\dots\ai{n})\rho_{\R,n}(A^{(n)})\ ,
\\
\lb{ae2}
(\ai{1}\dots\ai{n})\rho_{\R,n}(A^{(n)}) &=&
\rho_{\Rp,n}(A^{(n)})(\ai{1}\dots\ai{n})\rho_{\R,n}(A^{(n)})\ .
\ea
Finally, expressing (\ref{A}), (\ref{DA1}) for constant and
dynamical $q$-anti\-sym\-met\-ri\-zers in terms of the $\e$-tensors,
one transforms (\ref{ea2}), (\ref{ae2}) to the form (\ref{ea}), (\ref{ae}).
\eod
\vspace{2mm}

\noindent
{\bf Proposition 5.2} {\em
The element $\det (a)$ of the algebra ${\cal A}(\Rp, \R)$ commutes with
the generators $p_i$ and its commutation with the generators $a^i_{\a}$
is described by
\be
\lb{deta-a}
\det (a) \, a = K(p)\, a \det (a) \ ,
\ee
where the diagonal matrix $K(p)$ is given in (\ref{Kp}), (\ref{npi}).
}
\vspace{2mm}

\noindent
{\bf Proof.} Consider the permutation of $\det (a)$ with an
arbitrary function $h(p)$:
\ba
\nn
\lefteqn{
\det (a) h(p) =  \Edo\ai{1}\dots\ai{n}h(p)\eup /\raisebox{-2pt}{$\q{n}!$}
}&&
\\
\lb{transform}
&&=  \Edo(X_n\dots X_1)h(p)(X_n\dots X_1)^{-1}
\ai{1}\dots\ai{n}\eup /\raisebox{-2pt}{$\q{n}!$}
\ea
Since the only nonvanishing components of the tensor $\Edo$ are those
with pairwise  different indices and due to the diagonal
structure of the matrix $X$ one has
$$
\Edo X_n\dots X_1 = \Edo \det (X) = \Edo\ ,
$$
where in the last equality the unimodularity of $X$
(see (\ref{X})) is taken into account.
Now we can complete the transformation of (\ref{transform}):
$$\det (a) h(p) = \ldots $$
$$
=  h(p) \Edo(X_n\dots X_1)^{-1}
\ai{1}\dots\ai{n}\eup /\raisebox{-2pt}{$\q{n}!$} =
h(p) \det (a)\ .
$$
This proves commutativity of $\det (a)$ and $p_i$.

Consider now permutation of $\det (a)$ with the matrix $a$.
It is technically convenient to take $a$ living in the matrix
space with label $(n+1)$:
\ba
\nn
\lefteqn{
\det (a)\, \ai{n{+}1} =   \Edo \left\{\ai{1}\dots\ai{n}
\ai{n{+}1}\right\} \eup /\raisebox{-2pt}{$\q{n}!$}
}&&
\\
\nn
&&
= \left\{ \Edo\, \rho_{\Rp,n{+}1}(g_n\dots g_1)\right\}
\ai{1}\dots \ai{n{+}1}
\left\{\rho_{\R,n{+}1}^{-1}(g_n\dots g_1)\,\eup\right\}
/\raisebox{-2pt}{$\q{n}!$}
\\
\nn
&& =  \K{n{+}1}{1}(p)\left\{ X_1 \Ed2 X_1^{-1}\,\ai{1}\right\}
\left\{\ai{2}\dots\ai{n{+}1} \eu2\right\}
\iK{1}{n{+}1} /\raisebox{-2pt}{$\q{n}!$}
\\
\nn
&& = \K{n{+}1}{1}(p)\, \ai{1}\left\{\Ed2 \Eu2\right\}
\det (a)\,
\iK{1}{n{+}1} /\raisebox{-2pt}{$\q{n}!$}
\\
\nn
&& = \left( K(p) a K^{-1}\right)_{n{+}1} \det (a)\ .
\ea
The following formulae are
used in the course of the calculation: (\ref{deta})
and (\ref{intertwine})
in the first line, (\ref{DeR}) and (\ref{R11c})	 in the second
line, (\ref{ap}) and (\ref{ae})  in the third line, and (\ref{norm})
in passing to the last line.
For clarity we put into braces those expressions which are to be
transformed in the next step.

Finally, substituting $\id$ for the matrix $K$ (see (\ref{nkid}))
we obtain (\ref{deta-a}).
\eod
\vspace{2mm}

\noindent
{\bf Corollary 5.1} {\em
The element
\be
\lb{Central}
\Delta = \det (a) \, \prod_{i<j}  {\varphi_{ij}(p_{ij})\over f(p_{ij})}\ ,
\ee
where $f(p_{ij}) = \bq^{p_{ij}} + \q{p_{ij}}\b_{ij}$
and the functions $\varphi_{ij} $ are defined by the relations
\be
\lb{varphi}
\a_{ij}(p_{ij})	= {\varphi_{ij}(p_{ij}+1)\over \varphi_{ij}(p_{ij})}
\ee
belongs to the center of the algebra
${\cal A}(\Rp, \R)$.
The $SL(n)$-reduction in the algebra ${\cal A}(\Rp, \R)$ can be
performed by imposing the condition
$\Delta = 1$.
}
\vspace{2mm}

\noindent
{\bf Proof.}  We shall search for the central element in
${\cal A}(\Rp, \R)$	in the form $\Delta = U(p) \det (a)$,
where $U(p)$ is some function of $p_i$ which is to be fixed.
As follows from the Proposition 5.2 the element $\Delta$
commutes with $p_i$ and its commutativity with the generators
$a^i_\a$ imposes the following conditions on the function $U$
\be
\lb{cond}
X^i U(p) (X^i)^{-1}	= U(p) K^i_i(p)\ , \qquad i=1,\dots ,n\ .
\ee
Now using (\ref{0.8}), (\ref{xi}) and (\ref{npi})
it is straightforward to check that with the choice
(\ref{Central}), (\ref{varphi})	one satisfies conditions (\ref{cond}).
\eod
\vspace{2mm}

\noindent
{\bf Proposition 5.3}	{\em Let the algebra ${\cal A}(\Rp,\R)$ be completed
with the inverse determinant of $a$:
$(\det a)^{-1} \det(a) = \det (a) (\det a)^{-1} = 1$.
Then the left and right inverse of $a\;$ is given by
}
\be
\lb{ainv}
\ainv{1}{n{+}1} = {(-1)^{n-1}\over \q{n-1}!}\, (\det a)^{-1}\,
\Ed2\, \ai{2} \dots \ai{n} \eup\ .
\ee

\noindent
{\bf Proof.} We first check that the expression (\ref{ainv})
is a left inverse of $a\,$:
\ba
\nn
\lefteqn{
\ainv{1}{n+1}\, \ai{n+1} = {(-1)^{n-1}\over \q{n-1}!}\, (\det a)^{-1}
\left\{\Ed2 \ai{2}\dots \ai{n}\ai{n+1}\right\}\eup
}&&
\\
\nn
&& =
{(-1)^{n-1}\over \q{n-1}!}\, (\det a)^{-1}	\det (a)\, \ed2 \eup
\, =\, \N{1}{n+1} \, = \, \ID{1}{n+1}\ .
\ea
Here we have used successively Eqs.(\ref{ea}),
(\ref{N}) and (\ref{nkid}).

Checking that (\ref{ainv})
is also a right inverse is slightly more complicated:
\ba
\nn
\lefteqn{
\ai{1}\ainv{1}{n+1}	= {(-1)^{n-1}\over \q{n-1}!}
\left\{ \ai{1} (\det a)^{-1}\right\} \Ed2\, \ai{2}\dots \ai{n}\eup
}
\\
\nn
&&=
{(-1)^{n-1}\over \q{n-1}!}\, (\det a)^{-1} K_1(p)
\left\{ \ai{1} \Ed2 \right\} \ai{2}\dots \ai{n}\eup
\\
\nn
&&=
{(-1)^{n-1}\over \q{n-1}!}\, (\det a)^{-1} K_1(p)\, X_1\,
\Ed2 X_1^{-1}
\left\{ \ai{1}\dots	\ai{n}\eup\right\}
\\
\nn
&&=
(\det a)^{-1}  K_1(p)\, \N{1}{n+1}(p)\,\det (a) = \ID{1}{n+1}\ ,
\ea
where we have applied successively Eqs.(\ref{deta-a}),
(\ref{ap}), (\ref{ae}),
(\ref{Np}) and (\ref{DKN}).
\hfill\eod

\vspace{2mm}

The existence of inverse matrix $a^{-1}$ is needed in
many applications of the algebra ${\cal A}(\Rp , \R)$.
As an example of such application we shall construct a realization of a
reflection equation algebra ${\cal M}(\R)$
(for definition of this algebra see
e.g. \cite{Kul} and references therein)
in terms of the generators of
${\cal A}(\Rp, \R)$. We have to use here the following
general property of $SL(n)$-type dynamical $R$-matrices
(which has been noticed in \cite{AF} for the
$SL(2)$ case, see also \cite{BF}, \cite{By}):

\vspace{2mm}

\noindent
{\bf Proposition 5.4}	{\em
The dynamical matrix $\Rp$ (\ref{Rp}), (\ref{ab}),
(\ref{xi}) satisfies the equation
\be
\lb{mon1}
D_1 \, \Rp \, {D_2}^{-1} =
\Rp^{-1} \, \sigma_{12} \; ,
\ee
where the diagonal matrices $D$ and $\sigma$
\be
\lb{M3}
D^{i}_{j} \equiv q^{d_i} \d^{i}_{j} \;\; , \;\;\;
({\sigma}_{12})^{i_1 i_2}_{j_1 j_2} = \d^{i_1}_{j_1} \d^{i_2}_{j_2} \,
\sigma_{i_1 i_2} \; .
\ee
are fixed by (\ref{mon1}) as
\be
\lb{M6}
q^{d_i-d_j} = q^{- 2p_{ij}} \, \pi_{ij} \;\;\; (i \neq j) \;\, ,
\ee
\be
\lb{M7}
\sigma_{i j} = q^{2 \d_{ij}} \;
\ee
}
so that $d_i$ are functions of $p\,$.

\vspace{2mm}

\noindent
{\bf Proof.} First of all we note that from the Hecke condition
(\ref{dhecke}) (and (\ref{btt})) one can deduce
\be
\lb{mon3}
\Rp^{-1} =
\left( a_{i_1 i_2}(p) - (q-\bar q) \d_{i_1 i_2} \right)
\, \d^{i_1}_{j_2} \d^{i_2}_{j_1}\, - \,
b_{i_2 i_1}(p)\, \d^{i_1}_{j_1} \d^{i_2}_{j_2}\ ,
\ee
Substitution of (\ref{Rp}), (\ref{M3})
and (\ref{mon3}) into (\ref{mon1}) gives the
following equations for the parameters
$\sigma_{ij}$ and $d_i$
\be
\lb{mon4}
a_{ij} = \left( a_{ij} - (q - \bar q) \d_{ij} \right) \,
\sigma_{ji} \; ,
\ee
\be
\lb{mon5}
q^{d_i - d_j} \, b_{ij} =
- b_{ji} \, \sigma_{ij} \; .
\ee
Equation (\ref{mon4}) leads to (\ref{M7})
while (\ref{mon5}) is equivalent (in view of (\ref{is})) to (\ref{M6}).
\eod

\vspace{0.5cm}

Now we construct the matrix ${M^{\alpha}}_{\beta}$ which
is diagonalized
with the help of the matrix $a^{i}_{\alpha}$ and the spectrum of which
is defined by the matrix $D$ (\ref{M3}), (\ref{M6})
\be
\lb{M}
M = a^{-1} \, D \, a \;\;  .
\ee

It is clear that
$[D_1 , \, D_2] =0$ and therefore
the spectrum of the matrix $M$ gives a commutative set of elements.

\vspace{2mm}

\noindent
{\bf Proposition 5.5}	{\em
The elements of the matrix $M$ (\ref{M})
satisfy a reflection equation of the form
\be
\lb{M5}
M_2 \, \R^{-1} \, M_2 \, \R^{-1} =
\R^{-1} \, M_2 \, \R^{-1} \, M_2 \; ,
\ee
and thus provide a realization of a reflection equation
subalgebra ${\cal M}(\R)$ in ${\cal A}(\Rp, \R)$.
The matrix elements of $M$ satisfy the following exchange relations with
the generators of ${\cal A}(\Rp , \R)$:
\be
\lb{M4}
[ D_2 , \, M_1 ] =0 \;\; , \;\;\;
M_1 \, a_2 = q^{2/n} \,
a_2 \, \R^{-1} \, M_2 \, \R^{-1} \; .
\ee
}

\noindent
{\bf Proof.}
Using (\ref{ap}), one can bring the commutation relations of the
matrix $D$ with the elements $a^{i}_{\alpha}$
to the form
\be
\lb{M2}
a_1 \, D_2  = q^{-2/n} \, \sigma_{12} \,
D_2 \, a_1 \; ,
\ee
where the diagonal matrix $\sigma_{12}$ is given by (\ref{M3}), (\ref{M7}).
Eqs.(\ref{M}) and (\ref{M2}) imply
$[D_2 , \, M_1] =0$.
Then one proves (\ref{M5}) and the second relation in (\ref{M4}) by using
(\ref{M}), (\ref{aa}), (\ref{M2}) and (\ref{mon1})~.
\eod

\vspace{0.5cm}

\section{Application to the $SU(n)$ WZNW model}
\setcounter{equation}{0}
\renewcommand\theequation{\thesection.\arabic{equation}}

\medskip

As an application of quantum matrix algebras we
briefly describe here a typical problem of the two
dimensional conformal field theory in which such matrices arise (see
\cite{FHIPT} for more details).

Let $G$ be a connected compact Lie group and $g=g(t,x)$ be a map from the
cylinder ${\Bbb R}\times{\Bbb S}^1$ into $G$ which satisfies the
Wess--Zumino--Novikov--Witten (WZNW) equations of motion. The general
periodic solution $g(t,x) = g(t,x+2\pi ) $ of these equations
factorizes into a product of group valued chiral fields
\be
g^A_B(t,x)=u^A_\alpha (x-t){\bar u}^\alpha_B (x+t)\quad
({\rm classically,}\ \ g, u, {\bar u}\in G),
\lb{5.1}
\ee
each of which satisfies a twisted periodicity condition; in particular,
\be
u(x+2\pi )=u(x)M\,,\quad \left( M\in G\right)
\lb{5.2}
\ee
where $M$ is the monodromy.

Furthermore, the quantum chiral fields obey quadratic exchange relations
\cite{B, F1, F11, F2, AF, G, CG, FHT1, FHT2}
\be
\lb{5.3}
u(y)_2 u(x)_1 =\, u(x)_1 u(y)_2 R(x-y)\,\,
\Leftrightarrow \,\, P u(y)_2 u(x)_1 =\, u(x)_2 u(y)_1
\R (x-y)\,.
\ee
Here the matrix $R(x)$ is a solution of the
the quantum Yang--Baxter equation whose $x$-dependence
is given by a step function, while $\R (x)$ is the associated
braid operator:
\be
\lb{5.4}
\R (x) = \R \,\theta (x) + {\R}^{-1} \,\theta (-x)\,,
\quad \R (x) = P\, R(x)\, = {\R}^{\e (x)}
\ee
$\left(\e (x) = \theta (x) - \theta (-x)\,\right)\,.$

Since $\R$ enters Eq.(\ref{5.3}) in pair with $P$ it should be normalized
to have determinant ${\rm det}\, \R = {\rm det}\, P$. For $G=SU(n)$ this
implies the relation
\be
\lb{5.5}
{\R}_{i i+1} = {\bar q}^{1\over n}\,\rho (g_i )\
\left(\;{\rm for}\  g^2_i = \id + (q- {\bar q}) g_i\  \right)\
\Longrightarrow\
{\rm det}\, \R =
{\rm det}\, P = (-1)^{{n\choose 2}}
\ee
so that we have to renormalize $\R$ of (\ref{Jimbo}) by multiplying it by
${\bar q}^{1\over n}\,$. (The resulting $\R$ has eigenvalues
$q^{1-{1\over n}}$ and ${- \bar q}^{1+ {1\over n}}$ of multiplicities
${n+1\choose 2}$ and ${n\choose 2}$, respectively; thus the product
of all $n^2$ eigenvalues of $\R$ is indeed $(-1)^{{n\choose 2}}\,$.)

We expand, following \cite{FHT2} and \cite{BF}, $u(x)$ into a basis of
zero modes that diagonalizes the monodromy matrix $M$:
\be
\lb{5.6}
u^A_\a (x) = a^i_\a\, u^A_i (x,{p})\,,\quad
a\, M= D\, a\,,\quad D^i_j = q^{d_i} {\delta}^i_j\,.
\ee
Here $d_i = -2{p}_i -1/n + 1$,
$p = \{{p}_i\}$ are central elements of the reflection
equation algebra ${\cal M}(\R)\,$; in the quantum field
theoretic representation at hand they form a commuting set of operators
such that Eq.(\ref{0.7}) takes place.
The eigenvalues of the differences $p_{i i+1}\ ( = p_i - p_{i+1} )$
are natural numbers that can be identified with the extended weights,
${\lambda}_i + 1$ labeling the (finite dimensional) irreducible
representations of $SU(n)$ . The labels of the ${n\choose j}$
dimensional fundamental representation are given by
${\lambda}_i^{(j)} = \d^j_i\,,\ 1\le i,j\le n-1 $. Under these assumptions
Eq.(\ref{5.3}) implies exchange relations of the type
\be
\lb{5.7}
{\tilde R} (p)\, a_2\, a_1 = a_2\, a_1\,\R
\ee
where ${\tilde R} (p)\,$ obeys a QDYBE analogous to (\ref{0.3}).
Hence, the results displayed in Sections 3 and 4 can be applied with slight
modifications. (Since $\R$ and ${\tilde R} (p)\,$ enter (\ref{5.7})
homogeneously, the factor ${\bar q}^{1\over n}$ of
(\ref{5.5}) cancels in the two sides). Thus we can also
apply the results of Section 5 to the (chiral zero mode) quantum
matrix algebra ${\cal A}$ of the $SU(n)$ WZNW model. It should be noted
that in this case $q$ is a root of $-1$ associated with the level
$k\,\, {\widehat{su}}(n)$ Kac--Moody algebra:
\be
\lb{5.8}
q=e^{i{{\pi}\over h}}\,,\quad
[2] = q + \bq = 2 {\rm cos} {{\pi}\over h}\,,\quad
h=n+k \left( \ge n+1 \right)\,.
\ee

The eigenvalues $q^{d_i}$ of the diagonal matrix
$D$ can be expressed as differences of conformal dimensions. Indeed,
according to \cite{FHT2}, the chiral vertex operators $u_j(x, p)$ satisfy
\be
\lb{5.9}
u_j(x+2\pi , p) = u_j(x, p)\,
e^{2\pi i \left(\Delta_h({p})-\Delta_h({p} +v^{(j)}) \right)}\,,
\ee
where the matrices $v^{(j)}\,$ and $p\,$ are defined by
(\ref{0.1}) and (\ref{0.2}). Here the conformal dimensions are expressed
in terms of the $SU(n)$ Casimir operator,
\be
\lb{5.10}
2\, h\, \Delta_h (p) = C_2 (p) = {1\over n}
\sum\limits_{1\le i<k\le n} p^2_{ik}\, -\,{{n(n^2 -1)}\over 12}\,,
\ee
so that
\be
\lb{5.11}
d_j = C_2 (p) - C_2 (p+v^{(j)} )
= -2 ( p | v^{(j)} ) - |v^{(j)}|^2
= {1\over n} - 1 - 2 p_j\,.
\ee
(cf. (\ref{M6})).

An important consequence of
(\ref{5.7}) and (\ref{5.8})
is the existence of an ideal ${\cal I}_h$ of ${\cal A}$ generated by $n^2$
elements $(a^i_{\a})^h$ such that the factor algebra ${\cal A}/{\cal I}_h$
is finite dimensional \cite{FHIPT}~. This allows to define a finite
dimensional "Fock space representation" of ${\cal A}$ with a unique
vacuum vector $|vac>\,$ corresponding to trivial $su(n)$ weight
${\lambda}_i = 0\ (p_{i i+1} = 1\,,\ i=1, \ldots , n-1 )$ such that
\be
\lb{5.13}
a^i_{\a}\, |vac> = 0\quad {\rm for}\quad i> 1\,,\quad {\cal I}_h\, |vac>
= 0\,.
\ee


\section*{\bf Acknowledgements}

\medskip

This work started while three of us, L.K.H., O.V.O.
and I.T.T. were visiting the
Bogoliubov Laboratory of Theoretical Physics of the JINR, Dubna and was
completed while A.P.I., L.K.H. and I.T.T. were visiting the Department of
Physics of the University of Pisa, ICTP (Trieste) and the Erwin
Schr\"odinger Institute for Mathematical Physics, respectively.  We thank
all these Institutions for hospitality and support.  This work was
also supported in part by RFBR (grant 97-01-01041),
INTAS (grant 93-127-ext), by the exchange program
between INFN and JINR (Dubna) and by the Bulgarian National
Foundation for Scientific Research under contract F-404.

\vspace{0.5cm}


\section*{Appendix. Normalization of dynamical Levi-Civit\`a tensors}
\setcounter{equation}{0}
\def\theequation{A.\arabic{equation}}

\medskip

The definitions (\ref{Eup}), (\ref{Edo}) lead to the expression
\be
\lb{3.14c}
{\cal E}_{i_1 \dots i_n}(p) \, {\cal E}^{i_1 \dots i_n}(p) =
 \prod_{1 \leq a<b \leq n} \, \xi_{i_a i_b} \;
\ee
(there are no summations over the indices $i_k$), and the normalization
condition (\ref{norm}) for the dynamical ${\cal E}$-tensors
follows from
\\

\noindent
{\bf Proposition A~} {\em Let $\xi_{ij} = d - b_{ij}$ where $d$ is a
constant (comparing with (\ref{ab}), one gets $d = q$) and the elements
$b_{ij}$ satisfy (\ref{btt}), (\ref{btt2}). Then the following identity
holds:
\be
\lb{A.1}
I_k \equiv \sum_{S_k} \,
\prod_{1 \leq a < b \leq k} \, \xi_{i_a i_b} =
[k]_{d}! \; ,
\ee
where $k \leq n$,
$[k]_{d} = \frac{(d^{k} - (d-\lambda)^{k})}{\lambda} \ \ (\lambda = q-\bar q)$
and $S_k$ denotes all permutations of indeces $(i_1, \dots , i_k)$
($i_a \neq i_b$ for $a \neq b$).
Note that $[k]_{d} =[k]$
for $d=q$ (as it is needed in (\ref{norm})). }

\vspace{0.2cm}

\noindent
{\bf Proof~} We
shall proceed by induction. For $k=2$ we have
$$
I_2 = \xi_{i_1 i_2} + \xi_{i_2 i_1} = 2 \, d - \lambda =
[2]_{d} \; .
$$
Let (\ref{A.1}) is correct for some $k \; (1 < k < n)$,
then for $k+1$ we derive
$$
I_{k+1} = \sum_{S_{k+1}} \, \left[
( \prod_{l=1}^{k} \, \xi_{i_l i_{k+1}} ) \,
\prod_{1 \leq a < b \leq k} \, \xi_{i_a i_b} \right] =
$$
$$
= \sum_{r=1}^{k+1} \left[ ( \prod_{\stackrel{l=1}{l \neq r}}^{k+1}
\xi_{i_l i_r} ) \, \sum_{S_k}
\prod_{
\stackrel{a \neq r \neq b}{1 \leq a < b \leq k+1}
}
\, \xi_{i_a i_b} \right]
=  I_k \, \sum_{r=1}^{k+1}
 \prod_{\stackrel{l=1}{l \neq r}}^{k+1}
\xi_{i_l i_r}   \; .
$$
Therefore we should prove the identity
\be
\lb{A.2}
[k+1]_{d} =
\sum_{r=1}^{k+1}  \prod_{\stackrel{l=1}{l \neq r}}^{k+1}
\xi_{i_l i_r}  =
\sum_{r=1}^{k+1}  \prod_{\stackrel{l=1}{l \neq r}}^{k+1}
(d  -  b_{i_l i_r})  \; .
\ee
This identity follows from the relation
\be
\lb{A.3}
\sum_{r=1}^{m} \, \prod^{m}_{\stackrel{l=1}{l \neq r}} b_{i_l i_r} =
\lambda^{m-1} \;\; (m \leq k+1) \; ,
\ee
which can be obtained by induction. Indeed,
from Eqs.(\ref{btt}), (\ref{btt2})
we have for $m=2,3$
$$
b_{i_1 i_2} + b_{i_2 i_1} = \lambda \;\; , \;\;\;
b_{i_2 i_1} \, b_{i_3 i_1}  + b_{i_1 i_2} \, b_{i_3 i_2} +
b_{i_1 i_3} \, b_{i_2 i_3} = \lambda^{2} \; .
$$
Then we deduce
$$
\prod_{l=2}^{m} \, b_{i_l i_1} =
b_{i_2 i_1} \, \left(
\lambda^{m-2} - \sum_{r=3}^{m} \, (\lambda - b_{i_r i_1}) \,
\prod_{\stackrel{l=3}{l \neq r}} b_{i_l i_r} \right) =
$$
$$
= b_{i_2 i_1} \, \sum_{r=3}^{m} \, b_{i_r i_1} \,
\prod_{\stackrel{l=3}{l \neq r}} b_{i_l i_r}  =
\sum_{r=3}^{m} \, (\lambda^{2}
- b_{i_1 i_2} \, b_{i_r i_2}
- b_{i_1 i_r} \, b_{i_2 i_r} ) \,
\prod_{\stackrel{l=3}{l \neq r}} b_{i_l i_r}  =
$$
$$
= \lambda^{m-1} -
\sum_{r=2}^{m} \, \prod_{\stackrel{l=1}{l \neq r}} b_{i_l i_r}  \; ,
$$
which proves (\ref{A.3}). Expanding the right hand side of
(\ref{A.2}) in power series of $d$ and taking into account
(\ref{A.3}) we verify the relations (\ref{A.2}) and, thus,
complete the proof.	\eod

One can reformulate the statement of Proposition A
in more concise form (only in terms of elements $\xi_{ij}$)\\

\noindent
{\bf Proposition B~} {\em Let $\xi_{ij}$ satisfy
$$
\xi_{ij} + \xi_{ji} = [2] =
\xi_{ij} \, \xi_{jk} \, \xi_{ki} +
\xi_{ik} \, \xi_{kj} \, \xi_{ji} \;\; (i \neq j \neq k \neq i) \; .
$$
We rewrite these conditions as
\be
\lb{pro}
\sum_{r=1}^{k}  \prod_{\stackrel{l=1}{l \neq r}}^{k}
\xi_{i_l i_r}  =
[k] \;\; {\rm for} \;\; k=2,3 \; .
\ee
Then, equation (\ref{pro}) is also valid for
$4 \leq k \leq n$, and the following identity holds:
\be
\lb{pro1}
I_k \equiv \sum_{S_k} \,
\prod_{1 \leq a < b \leq k} \, \xi_{i_a i_b} =
[k] ! \; ,
\ee
where $S_k$ denotes all permutations of the indices
$(i_1, \dots , i_k)$ and $i_a \neq i_b$ for $a \neq b$.}

\vspace{0.2cm}

\noindent
{\bf Proof~} The proof is similar to that of Proposition A.	\eod
\vspace{0.2cm}

\noindent
{\bf Remark~} There are many other interesting relations among the elements
$b_{ij}$ (\ref{btt}), (\ref{btt2}) (as well as among $\xi_{ij}$). For
example, one can easily deduce the identity
$$
b_{i_1 i_2} \,
b_{i_2 i_3} \dots
b_{i_{k-1} \, i_k} \,
b_{i_{k} \, i_1} =
(-1)^{k} \, b_{i_1 i_k} \,
b_{i_k i_{k-1}} \dots
b_{i_{3} \, i_2} \,
b_{i_{2} \, i_1} \; ,
$$
which generalizes (\ref{btt2}) and follows from the relation
\be
\lb{is}
- \frac{b_{ji}(p)}{b_{ij}(p)} \, q^{2p_{ij}} = \pi_{ij} \; .
\ee
Note that we consider $\pi_{ij}$ as constants
which are independent of $p_{i}$.



\begin{thebibliography}{000}

\bibitem{FRT}
L.D. Faddeev, N. Yu. Reshetikhin and L.A. Takhtajan,
Algebra i Analiz {\bf 1:1} (1989) 178
(English translation: Leningrad Math. J. {\bf 1} (1990) 193).

\bibitem{GN}
J.-L. Gervais and A. Neveu, Nucl. Phys. {\bf B238} (1984) 125.

\bibitem{Skl}
E. Sklyanin,
Algebra i Analiz {\bf 6} (1994) 227 (English translation:
St. Petersburg Math. J. {\bf 6:2} (1995) 397-406), {\tt hep-th/9308060}.

\bibitem{BF}
A.G. Bytsko and L.D. Faddeev,
J. Math. Phys. {\bf 37} (1996) 6324, {\tt q-alg/9508022}.

\bibitem{AF}
A.Yu. Alekseev and L.D. Faddeev,
Commun. Math .Phys. {\bf 141} (1991) 413-422.

\bibitem{F1}
L.D. Faddeev, Commun. Math. Phys. {\bf 132} (1990) 131.

\bibitem{KR}
A.N. Kirillov and N.Yu. Reshetikhin,
{\em Representations of the algebra $U_q(sl(2)) ,\, q$-orthogonal
polynomials and invariants of links}, {\bf in:} Infinite Dimensional Lie
Algebras and Groups, Proc. 1988 CIRM Conference, Ed. V. Kac (World
Scientific, Singapore 1989), pp. 285-339.

\bibitem{BBB}
O. Babelon, D. Bernard and E. Billey,
Phys. Lett. {\bf B375} (1996) 89-97, {\tt q-alg/9511019};\\
O. Babelon, Commun. Math. Phys. {\bf 139 } (1991) 619-643.

\bibitem{W}
E. Witten,
Commun. Math. Phys. {\bf 92} (1984) 455.

\bibitem{F11}
A. Alekseev and S. Shatashvili, Commun. Math. Phys. {\bf 133}
(1990) 353.

\bibitem{F2}
L.D. Faddeev, Garg\`ese lectures 1991, {\bf in:}
New Symmetry Principles in Quantum Field
Theory, Ed. J. Fr\"ohlich et al., Plenum Press, NY (1992),
pp. 159--175.

\bibitem{Mai}
J.-M. Maillet,
Phys. Lett. {\bf B162} (1985) 137.

\bibitem{AT}
A.Yu. Alekseev and I.T. Todorov,
Nucl.Phys. {\bf B421} (1994) 413.

\bibitem{By}
A.G. Bytsko,
{\em Fusion of $q$-tensor operators: quasi-Hopf-algebraic point
of view}, {\tt q-alg/9609007}.

\bibitem{BS}
A.G. Bytsko and V. Schomerus,
{\em Vertex operators -- from a toy model to lattice algebras},
{\tt q-alg/9611010}.

\bibitem{Fe} G.Felder,
{\em Elliptic quantum groups,} {\bf in:} Proceedings
of the ICMP, Paris 1994. pp. 211-218. Ed. D. Iagolnitzer.
International Press, Boston, 1995, and
{\tt hep-th/9412207};
{\em Conformal field theory and integrable
systems associated with elliptic curves,} {\bf in:} Proceedings
of the ICM, Zurich, 1994, and
{\tt hep-th/9407154}.

\bibitem{ABB}
J.Avan, O.Babelon and E.Billey,
Comm. Math. Phys. {\bf 178} (1996) 281;\\
G. Arutyunov, L. Chehov and S. Frolov, {\em $R$-matrix quantization of
the elliptic Ruijsenaars--Schneider model}, {\tt q-alg/9612032}; \\
G. Felder and A.Varchenko, {\em
Elliptic
quantum groups and
Ruijsenaars models}, {\tt q-alg/9704005}.

\bibitem{I2}
A.P. Isaev,
J. Phys. A: Math. Gen. {\bf 29} (1996) 6903.

\bibitem{EV} P. Etingof and A. Varchenko,
{\em Solutions of the quantum dynamical Yang-Baxter
equation and dynamical quantum groups}, {\tt q-alg/9708015}.

\bibitem{FHT1}
P. Furlan, L.K. Hadjiivanov and I.T. Todorov,
{\em Canonical approach to the quantum WZNW model},
ICTP Trieste and ESI Vienna preprint IC/95/74, ESI 234 (1995).

\bibitem{FHT2}
P.Furlan, L.K.Hadjiivanov and I.T. Todorov,
Nucl. Phys. {\bf B474} (1996) 497, {\tt hep-th/9602101};\\
P.Furlan, L.K.Hadjiivanov and I.T. Todorov,
Int. J. Mod. Phys. {\bf A12} (1997) 23, {\tt hep-th/9610202}.

\bibitem{DT}
M. Dubois-Violette and I.T. Todorov,
Lett. Math. Phys. {\bf 42} (1997) 183, {\tt hep-th/9704069}.

\bibitem{G}
K. Gaw\c{e}dzki, Commun. Math. Phys. {\bf 139} (1991) 201;\\
F. Falceto and K. Gaw\c{e}dzki, J. Geom. Phys. {\bf 11} (1993)
251.

\bibitem{FHIPT}
P. Furlan, L.K. Hadjiivanov, A.P. Isaev, P.N. Pyatov and I.T. Todorov
(in preparation).

\bibitem{Gur} D.I. Gurevich,
{\em Algebraic aspects of the quantum Yang-Baxter equation},
Algebra i Analiz {\bf 2} (1990) 119-148
(English translation: Leningrad Math. J. {\bf 2}
(1991) 801-828).

\bibitem{M}
P. Martin,
{\em Potts Models and Related Problems in Statistical Mechanics,}
World Scientific, Singapore, 1991.

\bibitem{D} V.G. Drinfeld,
{\em Quantum groups,} {\bf in:}
Proceedings of the International Congress of Mathematicians,
Berkeley 1986, Academic Press (1986), Vol.1, pp. 798-820.

\bibitem{J}
M. Jimbo,
Lett.Math.Phys. {\bf 11} (1986) 247.

\bibitem{Drinfel} V.G.Drinfeld, Algebra i Analiz {\bf 1:6}
(1989) 141 (English translation: Leningrad Math. J. {\bf 1} (1990)
1419); \\
N.Yu.Reshetikhin, Lett. Math. Phys., {\bf 20} (1990) 331.

\bibitem{Man} Yu.I. Manin,
\it Quantum groups and noncommutative geometry, \rm Universit\'e de
Montr\'eal preprint CRM-1561 (1989);
\rm Comm. Math. Phys. \bf 122 \rm (1989) 163-175.

\bibitem{GIOPS}
P.N. Pyatov and  P.A. Saponov,
J. Phys. A: Math. Gen. {\bf 28}	(1995) 4415--4421;\\
D.I. Gurevich, P.N. Pyatov and P.A. Saponov,
Lett. in Math. Phys. {\bf 41} (1997) 255--264;\\
A.P. Isaev, O.V. Ogievetsky, P.N. Pyatov and P.A. Saponov,
{\em Characteristic polynomials and Newton identities for
quantum matrices},
Preprint CPT-97/P3471 (1997).

\bibitem{FeV} G.Felder and A.Varchenko,
Comm. Math. Phys. {\bf 181} (1996) 741.

\bibitem{CrGe}
E. Cremmer and J.-L. Gervais,
Comm. Math. Phys. {\bf 134} (1990) 619;\\
A. Bilal and J.-L. Gervais,
Nucl. Phys. {\bf B318} (1989) 579.

\bibitem{Kul} P.P. Kulish and R. Sasaki,
Prog.Theor.Phys. {\bf 89} \rm (1993) 741.

\bibitem{B}
O. Babelon, Phys. Lett. {\bf B215} (1988) 523;\\
B. Blok, Phys. Lett. {\bf B233} (1989) 359.

\bibitem{CG}
M. Chu, P. Goddard, I. Halliday, D. Olive and A. Schwimmer,
Phys. Lett. {\bf B266} (1991) 71;\\
M. Chu and P. Goddard,
Nucl. Phys. {\bf 445} (1995) 145.








\end{thebibliography}
\end{document}